\documentclass[12pt]{article}
\usepackage{amsmath, amssymb, amsfonts, amsthm, amscd, vmargin}
\usepackage[latin1]{inputenc}
\usepackage[all]{xy}
\usepackage{graphicx}
\usepackage{enumerate}

\setmarginsrb{2.4cm}{4.5cm}{2.4cm}{3.3cm}{0cm}{0cm}{0cm}{1.5cm}

\theoremstyle{plain}
\newtheorem{teo}{Theorem}[section]
\newtheorem{lem}[teo]{Lemma}
\newtheorem{prop}[teo]{Proposition}
\newtheorem{cor}[teo]{Corollary}

\theoremstyle{definition}
\newtheorem{defi}[teo]{Definition}

\newtheorem{rem}[teo]{Remark}

\numberwithin{equation}{teo}
\newcommand{\Aut}{\operatorname{Aut}}

\newcommand{\Lie}{\operatorname{Lie}}

\renewcommand{\dim}{\operatorname{dim}}

\renewcommand{\Im}{\operatorname{Im}}

\newcommand{\Cbb}{{\mathbb C}}

\newcommand{\Zbb}{{\mathbb Z}}
\newcommand{\Rbb}{{\mathbb R}}
\newcommand{\Pbb}{{\mathbb P}}

\newcommand{\Obb}{{\mathbb O}}

\newcommand{\lra}{\longrightarrow}
\newcommand{\lmt}{\longmapsto}

\begin{document}

\title{On some smooth projective two-orbits varieties with Picard number~1}

\author{Boris Pasquier}
\maketitle

\begin{abstract} We classify all smooth projective horospherical varieties with Picard
number~1. We prove that the automorphism group of any such variety
$X$ acts with at most two orbits and that this group still acts with
only two orbits on $X$ blown up at the closed orbit. We characterize
all smooth projective two-orbits varieties with Picard number~1 that
satisfy this latter property.
\end{abstract}

\textbf{Mathematics Subject Classification.} 14J45 14L30 14M17\\

\textbf{Keywords.} Two-orbits varieties, horospherical varieties,
spherical varieties.

\section*{Introduction}

Horospherical varieties are complex normal algebraic varieties where
a connected complex reductive algebraic group acts with an open
orbit isomorphic to a torus bundle over a flag variety. The
dimension of the torus is called the rank of the variety. Toric
varieties and flag varieties are the first examples of horospherical
varieties (see \cite{Pa06} for more examples and background).

It is well known that the only smooth projective toric varieties
with Picard number~1 are the projective spaces.  This is not the
case for horospherical varieties: for example any flag variety $G/P$
with $P$ a maximal parabolic subgroup of $G$ is smooth, projective
and horospherical with Picard number~1.

Moreover, smooth projective horospherical varieties with Picard
number~1 are not necessarily homogeneous. For example, let $\omega$
be a skew-form of maximal rank on $\Cbb^{2m+1}$. For
$i\in\{1,\ldots,m\}$, define the odd symplectic grassmannian
$\operatorname{Gr}_\omega(i,2m+1)$ as the variety of $i$-dimensional
$\omega$-isotropic subspaces of $\Cbb^{2m+1}$. Odd symplectic
grassmannians are horospherical varieties (see Proposition
\ref{oddgrass}) and, for $i\neq m$ they have two orbits under the
action of their automorphism group which is a connected
non-reductive linear algebraic group (see \cite{Mi05} for more
details).

Our focus on smooth horospherical varieties with Picard number~1,
that gives interesting examples of Fano varieties with Picard
number~1, is also motivated by the main result of \cite{Pa06}, where
Fano horospherical varieties are classified in terms of rational
polytopes. Indeed in \cite[Th.0.1]{Pa06}, the degree ({\it i.e.} the
self-intersection number of the anticanonical bundle) of smooth Fano
horospherical varieties is bounded. Two different bounds are
obtained in the case of Picard number~1 and in the case of higher
Picard number.

In Section \ref{section1}, we classify all smooth projective
horospherical varieties with Picard number~1.  More precisely, we prove the following result.

\begin{teo}\label{classifintro}
Let $G$ be a connected reductive algebraic group. Let $X$ be a
smooth projective horospherical $G$-variety with Picard number~1.

Then we have the following alternative:
\begin{enumerate}[(i)]
\item $X$ is homogeneous, or
\item $X$ is horospherical of rank~1. Its automorphism group is a connected non-reductive linear
algebraic group, acting with exactly two orbits.
\end{enumerate}
Moreover in the second case, $X$ is uniquely determined by its two
closed $G$-orbits $Y$ and $Z$, isomorphic to $G/P_Y$ and $G/P_Z$
respectively; and $(G,P_Y,P_Z)$ is one of the triples of the
following list.
\begin{enumerate}
\item $(B_m,P(\omega_{m-1}),P(\omega_m))$ with $m\geq 3$
\item $(B_3,P(\omega_1),P(\omega_3))$
\item $(C_m,P(\omega_i),P(\omega_{i+1}))$ with $m\geq 2$ and
$i\in\{1,\ldots,m-1\}$
\item $(F_4,P(\omega_2),P(\omega_3))$
\item $(G_2,P(\omega_2),P(\omega_1))$
\end{enumerate}
 Here we denote by $P(\omega_i)$ the maximal
parabolic subgroup of $G$ corresponding to the dominant weight
$\omega_i$ with the notation of Bourbaki \cite{Bo75}.
\end{teo}

Remark that Case 3 of Theorem \ref{classifintro} corresponds to odd
symplectic grassmannians. It would be natural to investigate other
complete smooth spherical varieties with Picard number~1 (A normal
variety is spherical if it admits a dense orbit of a Borel subgroup,
for example horospherical
varieties and symmetric varieties are spherical). A classification has been recently given in the special case of projective symmetric varieties by A~Ruzzi \cite{Ru07}.\\

In the second part of this paper, we focus on another special
feature of the non-homogeneous varieties classified by Theorem
\ref{classifintro}: the fact that they have two orbits even when
they are blown up at their closed orbit. Two-orbits varieties ({\it
i.e.} normal varieties where a linear algebraic group acts with two
orbits) have already be studied by D.~Akhiezer and S.~Cupit-Foutou.
In \cite{Ah83}, D.~Akhiezer classified those whose closed orbit is
of codimension~1 and proved in particular that they are
horospherical when the group is not semi-simple. In \cite{CF03},
S.~Cupit-Foutou classified two-orbits varieties when the group is
semi-simple, and she also proved that they are spherical. In section
\ref{section2}, we define two smooth projective two-orbits varieties
$\mathbf{X_1}$ and $\mathbf{X_2}$ with Picard number one (see
Definitions \ref{defiX1} and \ref{defiX2}) and we prove the
following:

\begin{teo}\label{caracterisation2}
Let $X$ be a smooth projective variety with Picard number~1 and put
$G:=\Aut^0(X)$.

Assume that $X$ has two orbits under the action of
$G$ and denote by $Z$ the closed orbit.
Then the codimension of $Z$ is at least $2$.

Assume furthermore that the blow-up of $Z$ in $X$ still has two
orbits under the action of $G$. Then, one of the following happens:
\begin{itemize}
\item $G$ is not semi-simple and $X$ is one of
the two-orbits varieties classified by Theorem \ref{classifintro};
\item $G=F_4$ and $X=\mathbf{X_1}$;
\item $G=G_2\times\operatorname{PSL(2)}$ and $X=\mathbf{X_2}$.
\end{itemize}
\end{teo}

The varieties in Theorem \ref{caracterisation2} are spherical of
rank one \cite{Br88}. Remark also that odd symplectic grassmannians
have been studied in detail by I.A.~Mihai in \cite{Mi05}. In
particular, he proved that an odd symplectic grassmannian is a
linear section of a grassmannian \cite[Prop.2.3.15]{Mi05}. It could
be interesting to obtain a similar description also for the
varieties of Theorem
\ref{caracterisation2}.\\

The paper is organized as follows.

In Section \ref{1}, we recall some results on horospherical
homogeneous spaces and horospherical varieties, which  we will use
throughout Section \ref{section1}. In particular we briefly
summarize the Luna-Vust theory \cite{LV83} in the case of
horospherical homogeneous spaces.

In Section \ref{2}, we prove that any horospherical
homogeneous space admits at most one smooth equivariant
compactification with Picard number~1. Then we give the list of
horospherical homogeneous spaces that admit a smooth
compactification not isomorphic to a projective space and with
Picard number~1. We obtain a list of 8 cases (Theorem
\ref{classif}).

In Section \ref{3}, we prove that in 3 of these cases, the smooth
compactification is homogeneous (under the action of a larger
group).

In Section  \ref{4}, we study the 5 remaining cases (they are listed
in Theorem \ref{classifintro}). We compute the automorphism group of
the corresponding smooth compactification with Picard number~1. We
prove that this variety has two orbits under the action of its
automorphism group and that the latter is connected and not
reductive.

In Section \ref{5}, we prove Theorem \ref{caracterisation2} in the
case where the automorphism group is not semi-simple. This gives
another characterization of the varieties obtained in Section
\ref{4}.\\

The aim of Section \ref{section2} is to prove Theorem
\ref{caracterisation2} when $G$ is semi-simple.
\begin{defi}\label{etoile}
A projective $G$-variety $X$ satisfies (*), if it is smooth with
Picard number~1, has two orbits under the action of $G$ such that
its closed orbit $Z$ has codimension at least~$2$, and the
blowing-up of $X$ along $Z$ has also two orbits under the action of
$G$.
\end{defi}
In Section \ref{genres}, we prove the first part of Theorem
\ref{caracterisation2} and we reveal two general cases.

In Sections \ref{RPH} and \ref{casiiiv}, we study these two cases
respectively. First, we reformulate part of the classification of
two-orbits varieties with closed orbit of codimension one due to
D.~N.~Akhiezer, in order to give a complete and precise list of
possible cases. Then we study separately all these possible cases.
We prove that the two varieties $\mathbf{X_1}$ and $\mathbf{X_2}$
satisfy (*) and are non-homogeneous, and that in all other cases,
the varieties satisfying (*) are homogeneous.

\section{Smooth projective horospherical varieties with Picard number~1}\label{section1}

\subsection{Notation}\label{1}

Let $G$ be a reductive and connected algebraic group over $\Cbb$, let
$B$ be a Borel subgroup of $G$, let $T$ be a maximal torus of $B$ and let $U$ be the
unipotent radical of $B$. Denote by $C$ the center of $G$ and by
$G'$ the semi-simple part of $G$ (so that $G=C.G'$). Denote by
$S$ the set of simple roots of $(G,B,T)$, and by $\Lambda$
(respectively $\Lambda^+$) the group of characters of $B$
(respectively the set of dominant characters). Denote by $W$ the
Weyl group of $(G,T)$ and, when $I\subset S$,  denote by $W_I$ the
subgroup of $W$ generated by the reflections associated to the
simple roots of $I$. If $\alpha$ is a simple root, we denote by
$\check\alpha$ its coroot,  and by $\omega_\alpha$ the  fundamental
weight corresponding to $\alpha$ (when the roots are
$\alpha_1,\dots\alpha_n$, we will write $\omega_i$ instead of
$\omega_{\alpha_i}$). Denote by $P(\omega_\alpha)$ the maximal
parabolic subgroup containing $B$ such that $\omega_\alpha$ is a
character of $P(\omega_\alpha)$. Let $\Gamma$ be the Dynkin diagram
of $G$. When $I\subset S$, we denote by $\Gamma_I$ the full subgraph
of $\Gamma$ with vertices the elements of $I$. For
$\lambda\in\Lambda^+$, we denote by $V(\lambda)$ the irreducible
$G$-module of highest weight $\lambda$ and by $v_\lambda$ a highest
weight vector of $V(\lambda)$. If $G$ is simple, we index the simple
roots as in
\cite{Bo75}.\\

A closed subgroup $H$ of $G$ is said to be {\it horospherical} if it
contains the unipotent radical of a Borel subgroup of $G$. In that
case we also say that the homogeneous space $G/H$ is {\it
horospherical}. Up to conjugation, one can assume that $H$ contains
$U$.  Denote by $P$ the normalizer $N_G(H)$ of $H$ in $G$. Then $P$
is a parabolic subgroup of $G$ such that $P/H$ is a torus. Thus
$G/H$ is a torus bundle over the flag variety $G/P$. The dimension
of the torus is called the {\it rank} of $G/H$ and denoted by $n$.

A normal variety $X$ with an action of $G$ is said to be a {\it
horospherical variety} if $G$ has an open orbit isomorphic to $G/H$
for some horospherical subgroup $H$. In that case, $X$ is also said
to be a $G/H$-embedding. The classification of $G/H$-embeddings (due
to D.~Luna et Th.~Vust \cite{LV83} in the more general situation of
spherical homogeneous spaces) is detailed in
\cite[Chap.1]{Pa06}.\\

Let us summarize here the principal points of this theory.  Let
$G/H$ be a fixed horospherical homogeneous space of rank $n$. This
defines a set  of simple roots $$I:=\{\alpha\in S\mid \omega_\alpha
\mbox{ is not a character of }P\}$$ where $P$ is the unique
parabolic subgroup associated to $H$ as above.  We also introduce a
lattice $M$ of rank $n$ as the sublattice of $\Lambda$ consisting of
all characters $\chi$ of $P$ such that the restriction of $\chi$ to
$H$ is trivial. Denote by $N$ the dual lattice to $M$.

In this paper, we call {\it colors} the elements of $S\backslash I$.
For any color $\alpha$, we denote by $\check\alpha_M$ the element of
$N$ defined as the restriction to $M$ of the coroot
$\check\alpha:\Lambda\lra\Zbb$. The point $\check\alpha_M$ is called
the {\it image of the color} $\alpha$. See \cite[Chap.1]{Pa06} to
understand the link between colors and the geometry of $G/H$.

\begin{defi}
A {\it colored cone} of $N_\Rbb:=N\otimes_\Zbb\Rbb$ is an ordered pair $(\mathcal{C}, \mathcal{F})$
 where $\mathcal{C}$ is a convex cone of $N_\Rbb$ and $\mathcal{F}$ is a set of colors (called the set of colors of the colored cone), such that\\
(i) $\mathcal{C}$ is generated by finitely many elements of $N$ and contains the image of the colors of $\mathcal{F}$,\\
(ii) $\mathcal{C}$ does not contain any line and
the image of any color of $\mathcal{F}$ is not zero.
\end{defi}
One defines a {\it colored fan} as a set of colored cones such that
any two of them intersect in a common colored face (see
\cite[def.1.14]{Pa06} for the precise definition).

Then $G/H$-embeddings are classified in terms of colored fans.
Define a {\it simple} $G/H$-embedding of $X$ as one containing a
unique closed $G$-orbit. Let $X$ be a $G/H$-embedding and
$\mathbb{F}$ its colored fan. Then $X$ is covered by its simple
subembeddings, and each of them corresponds to a colored cone of the
colored fan of $X$. (See \cite{Br97b} or \cite{Kn91} for the general
theory of spherical embeddings.)

\subsection{Classification of smooth projective embeddings with Picard number~1}\label{2}

The Picard number $\rho_X$ of a smooth projective $G/H$-embedding
$X$ satisfies
$$\rho_X=r_X+\sharp(S\backslash I)-\sharp(\mathcal{D}_X)$$
where $\mathcal{D}_X$ denotes the set of simple roots in
$S\backslash I$ which correspond to colors of $\mathbb{F}$ and $r_X$
is the number of rays of the colored fan of $X$ minus the rank $n$
\cite[(4.5.1)]{Pa06}. Since $X$ is projective, its colored fan is
complete ({\it i.e.} it covers $N_\Rbb$) and hence $r_X\geq 1$.
Moreover $\mathcal{D}_X\subset S\backslash I$, so $\rho_X=1$ if and
only if $r_X=1$ and $\mathcal{D}_X=S\backslash I$. In particular the
colored fan of $X$ has exactly $n+1$ rays.

\begin{lem}\label{unique}
Let $G/H$ be a horospherical homogeneous space.  Up to isomorphism of varieties, there exists at most one smooth projective $G/H$-embedding with Picard number~1.
\end{lem}
\begin{proof}
Let $X$ and $X'$ be two smooth projective $G/H$-embeddings with
respective colored fans $\mathbb{F}$ and $\mathbb{F'}$ and both with
Picard number~1.   Denote by $e_1,\ldots,e_{n+1}$ the primitive
elements of the $n+1$ rays of $\mathbb{F}$. By the smoothness
criterion of \cite[Chap.2]{Pa06}, $(e_1,\ldots,e_n)$ is a basis of
$N$, $e_{n+1}=-e_1-\cdots -e_n$ and the images in $N$ of the colors
are distinct and contained in  $\{e_1,\ldots,e_{n+1}\}$. The same
happens for $\mathbb{F'}$.  Then there exists an automorphism $\phi$
of the lattice $N$ which stabilizes the image of each color and
satisfies $\mathbb{F}=\phi(\mathbb{F'})$. Thus the varieties $X$ and
$X'$ are isomorphic  \cite[Prop. 3.10]{Pa06}.
\end{proof}

 If it exists, we denote by $X^1$ the unique smooth projective $G/H$-embedding with Picard number~1 and we say that $G/H$ is ``special''.

\begin{rem}
By the preceding proof, we have $\sharp(\mathcal{D}_X)\leq n+1$.
\end{rem}
\subsubsection{Projective space}

We first give a necessary condition for the embedding $X^1$ of a
special homogeneous space not to be isomorphic to a projective
space. In particular we must have $n=1$, so that $X^1$ has three
orbits under the action of $G$: two closed ones and $G/H$.

\begin{teo}\label{thesproj}
Let $G/H$ be a ``special'' homogeneous space. Then $X^1$ is
isomorphic to a projective space in the following cases:
\begin{description}
\item (i) $\sharp(\mathcal{D}_{X^1})\leq n$,
\item (ii) $n\geq 2$,
\item (iii) $n=1$, $\sharp(\mathcal{D}_{X^1})=2$ and the two simple roots of $\mathcal{D}_{X^1}$ are not in the same connected component of the Dynkin diagram $\Gamma$.
\end{description}
\end{teo}

\begin{proof}
(i) In that case, there exists a maximal colored cone of the colored
fan of $X$ which contains all colors. Then the corresponding simple
$G/H$-embedding of $X^1$, whose closed orbit is a point
\cite[Lem.2.8]{Pa06}, is affine \cite[th.3.1]{Kn91} and smooth.
 So it is necessarily a horospherical $G$-module $V$ \cite[Lem.2.10]{Pa06}. Thus $\Pbb(\Cbb\oplus V)$ is a smooth projective $G/H$-embedding with Picard number~1.
  Then by Lemma \ref{unique}, $X^1$ is isomorphic to $\Pbb(\Cbb\oplus V)$.\\

(ii) We may assume that $\sharp(\mathcal{D}_{X^1})=n+1$. Denote by
$\alpha_1,\ldots,\alpha_{n+1}$ the elements of $S\backslash I$ and
by $\Gamma_i$ the Dynkin diagram $\Gamma_{S\backslash\{\alpha_i\}}$.
The smoothness criterion of horospherical varieties
\cite[Chap.2]{Pa06} applied to $X^1$ tells us two things.

Firstly, for all $i\in\{1,\ldots,n+1\}$ and for all $j\neq i$,
$\alpha_j$ is a simple end ("simple" means not adjacent to a double
edge) of a connected component $\Gamma_i^j$ of $\Gamma_i$ of type
$A_m$ or $C_m$. Moreover the $\Gamma_i^j$ are distinct, in other
words, each connected component of $\Gamma_i$ has at most one vertex
among the $(\alpha_i)_{i\in\{1,\ldots,n+1\}}$.

Secondly, $(\check\alpha_{iM})_{i\in\{1,\ldots,n\}}$ is a basis of
$N$ and $\check\alpha_{(n+1)M}=-\check\alpha_{1M}-\cdots
-\check\alpha_{nM}$. Thus a basis of  $M$ (dual of $N$) is of the
form $$(\omega_i-\omega_{n+1}+\chi_i)_{i\in\{1,\ldots,n\}}$$ where
$\chi_i$ is a character of the center $C$ of $G$,  for all
$i\in\{1,\ldots,n\}$.

Let us prove that a connected component of $\Gamma$ contains at
most one vertex among the $(\alpha_i)_{i\in\{1,\ldots,n+1\}}$.
Suppose the contrary:  there exist $i,j\in\{1,\ldots,n+1\}$, $i\neq
j$ such that $\alpha_i$ and $\alpha_j$ are vertices of a connected
component of $\Gamma$. One can choose $i$ and $j$ such that there is
no vertex among the $(\alpha_k)_{k\in\{1,\ldots,n+1\}}$ between
$\alpha_i$ and $\alpha_j$. Since $n\geq 2$, there exists an integer
$k\in\{1,\ldots,n+1\}$ different from $i$ and $j$. Then we observe
that $\Gamma_k$ does not satisfy the condition that each of its
connected component has at most one vertex among the
$(\alpha_i)_{i\in\{1,\ldots,n+1\}}$ (because
$\Gamma_k^i=\Gamma_k^j$).

Thus we have proved that
\begin{equation}\label{eq1}
\Gamma=\bigsqcup_{j=0}^{n+1}\Gamma^j
\end{equation} such that for all $j\in\{1,\ldots,n+1\}$, $\Gamma^j$ is a connected component of $\Gamma$ of type $A_m$ or $C_m$ in which $\alpha_j$ is a simple end.

For all $\lambda\in\Lambda^+$, denote by by $V(\lambda)$ the simple $G$-module of weight $\lambda$.
Then Equation \ref{eq1} tells us that the projective space $$\Pbb(V(\omega_{n+1})\oplus V(\omega_1+\chi_1)\oplus\cdots\oplus V(\omega_n+\chi_n))$$
 is a smooth projective $G/H$-embedding with Picard number~1. Thus $X^1$ is isomorphic to this projective space.\\

(iii)  As in case (ii), one checks that $X^1$ is isomorphic to
$\Pbb(V(\omega_2)\oplus V(\omega_1+\chi_1))$ for some character
$\chi_1$ of $C$.
\end{proof}

\subsubsection{When $X^1$ is not isomorphic to a projective space}\label{2.2}

According to Theorem \ref{thesproj} we have to consider the case
where the rank of $G/H$ is~1 and where there are two colors
corresponding to  simple roots $\alpha$ and $\beta$  in the same
connected component of $\Gamma$. As we have seen in the proof of
Theorem \ref{thesproj}, the lattice $M$ (here of rank~1) is
generated by $\omega_\alpha-\omega_\beta+\chi$ where $\chi$ is a
character of the center $C$ of $G$. Moreover, $H$ is the kernel of
the character $\omega_\alpha-\omega_\beta+\chi:\,
P(\omega_\alpha)\cap P(\omega_\beta)\lra\Cbb^*$.

We may further reduce to the case where $G$ is semi-simple (recall
that $G'$ denotes the semi-simple part of $G$).

\begin{prop}
Let  $H'=G'\cap H$. Then $G/H$ is isomorphic to $G'/H'$.
\end{prop}

\begin{proof}
We are going to prove that $G/H$ and $G'/H'$ are both isomorphic to
a horospherical homogeneous space under $(G'\times \Cbb^*)$. In fact
$G/H$ is isomorphic to $(G'\times P/H)/\tilde{H}$ \cite[Proof of
Prop.3.10]{Pa06}, where $$\tilde{H}=\{(g,pH)\in G'\times P/H\mid
gp\in H\}.$$  Similarly, $G'/H'$ is isomorphic to $(G'\times
P'/H')/\tilde{H'}$ where $P'=P\cap G'$ and $\tilde{H'}$ defined as
the same way as $\tilde{H}$. Moreover the morphisms
$$\begin{array}{ccccccc}
P/H & \lra & \Cbb^*  & \mbox{ and } &  P'/H' & \lra & \Cbb^*\\
pH & \lmt & (\omega_\alpha-\omega_\beta+\chi)(p) &\mbox{~~~~~} & p'H' & \lmt & (\omega_\alpha-\omega_\beta)(p')
\end{array}$$
are isomorphisms. Then  \begin{eqnarray*} \tilde{H} & = &
\{(p',c)\in P'\times\Cbb^*\mid
(\omega_\alpha-\omega_\beta+\chi)(p')=c^{-1}\}\\
 & = & \{(p',c)\in P'\times\Cbb^*\mid
(\omega_\alpha-\omega_\beta)(p')=c^{-1}\}\\
& = & \tilde{H'}. \end{eqnarray*}
 This completes the proof.
\end{proof}

\begin{rem}\label{remP/H}
In fact $P/H\simeq\Cbb^*$ acts on $G/H$ by right multiplication, so
it acts on the $\Cbb^*$-bundle $G/H\lra G/P$ by multiplication on
fibers. Moreover, this action extends to $X^1$ (where $\Cbb^*$ acts
trivially on the two closed $G$-orbits).
\end{rem}

So we may assume that $G$ is semi-simple. Let $G_1,\ldots,G_k$ the
simple normal subgroups of $G$, so that $G$ is the quotient of the
product $G_1\times\cdots\times G_k$ by a central finite group $C_0$.
We can suppose that $C_0$ is trivial, because $G/H\simeq
\tilde{G}/\tilde{H}$ where $\tilde{G}=G_1\times\cdots\times G_k$ and
$\tilde{H}$ is the preimage of $H$ in $\tilde{G}$. If $\alpha$ and
$\beta$ are simple roots of the connected component corresponding to
$G_i$, denote by $H_i$ is the kernel of the character $\omega_\alpha
-\omega_\beta$ of the parabolic subgroup $P(\omega_\alpha)\cap
P(\omega_\beta)$ of $G_i$. Then
$$H=G_1\times\cdots\times G_{i-1}\times H_i\times
G_{i+1}\times\cdots\times G_k$$ and  $G/H=G_i/H_i$.

So from now on, without loss of generality, we suppose that $G$ is
simple.

\begin{teo}\label{classif}
With the assumptions above,  $G/H$ is "special" if and only if $(\Gamma,\alpha,\beta)$ appears in  the following list (up to exchanging $\alpha$ and $\beta$).
\begin{enumerate}
\item $(A_m,\alpha_1,\alpha_m)$, with $m\geq 2$;

\item $(A_m,\alpha_i,\alpha_{i+1})$, with $m\geq 3$ and $i\in\{1,\ldots,m-1\}$

\item $(B_m,\alpha_{m-1},\alpha_m)$, with $m\geq 3$;

\item $(B_3,\alpha_1,\alpha_3)$

\item $(C_m,\alpha_{i+1},\alpha_i)$ with $m\geq 2$ and $i\in\{1,\ldots,m-1\}$

\item $(D_m,\alpha_{m-1},\alpha_m)$, with $m\geq 4$;

\item $(F_4,\alpha_2,\alpha_3)$

\item $(G_2,\alpha_2,\alpha_1)$
\end{enumerate}
\end{teo}

\begin{proof}
The Dynkin diagrams $\Gamma_{S\backslash\{\alpha\}}$ and
$\Gamma_{S\backslash\{\beta\}}$ are respectively of type $A_m$ or
$C_m$ by the smoothness criterion \cite[Chap.2]{Pa06}. And for the
same reason, $\alpha$ and $\beta$ are simple ends of
$\Gamma_{S\backslash\{\beta\}}$ and $\Gamma_{S\backslash\{\alpha\}}$
respectively.

Suppose $\Gamma$ is of type $A_m$. If $\alpha$  equals $\alpha_1$
then, looking at $\Gamma_{S\backslash\{\alpha\}}$, we remark that
$\beta$ must be $\alpha_2$ or $\alpha_m$. So we are in Case 1 or 2.
If $\alpha$ equals $\alpha_m$ the argument is similar. Now if
$\alpha$ is not an end of $\Gamma$, in other words if
$\alpha=\alpha_i$ for some $i\in\{2,\ldots,m-1\}$ then, looking at
$\Gamma_{S\backslash\{\alpha\}}$, we see that $\beta$ can be
$\alpha_1$, $\alpha_{i-1}$, $\alpha_{i+1}$ or $\alpha_m$. The cases
where $\beta$ equals $\alpha_1$ or $\alpha_m$ are already done and
the case where $\beta$ equals $\alpha_{i-1}$ or $\alpha_{i+1}$ is
Case 2.

The study of the remaining cases is analogous and left to the
reader.
\end{proof}

In the next two sections we are going to study the variety $X^1$ for
each case of this theorem.  In particular we will see that $X^1$ is
never isomorphic to a projective space.

\subsection{Homogeneous varieties}\label{3}

In this section, with the notation of Section \ref{2.2}, we are going to prove that $X^1$ is homogeneous in Cases 1, 2 and 6.\\

In all cases (1 to 8), there are exactly 4 projective
$G/H$-embeddings and they are all smooth; they correspond to the 4
colored fans consisting of the two half-lines of $\Rbb$, without
color, with one of the two colors and with the two colors,
respectively (see \cite[Ex.1.19]{Pa06} for a similar example).

Let us realize $X^1$ in a projective space as follows.  The
homogeneous space $G/H$ is isomorphic to the orbit of the point
$[v_{\omega_\beta}+v_{\omega_\alpha}]$ in
$\Pbb(V(\omega_\beta)\oplus V(\omega_\alpha))$, where
$v_{\omega_\alpha}$ and $v_{\omega_\beta}$ are highest weight
vectors of $V(\omega_\alpha)$ and $V(\omega_\beta)$ respectively.
Then $X^1$ is the closure of this orbit in
$\Pbb(V(\omega_\beta)\oplus V(\omega_\alpha))$, because both have
the same colored cone ({\it i.e.} that with two colors)\footnote{See
\cite[Chap.1]{Pa06} for the construction of the colored fan of a
$G/H$- embedding.}.

We will describe the other $G/H$-embeddings in the proof of Lemma
\ref{lemaut}.

\begin{prop}
In Case 1, $X^1$ is isomorphic to the quadric $Q^{2m}=\operatorname{SO}_{2m+2}/P(\omega_1)$.
\end{prop}

\begin{proof}
Here, the fundamental $G$-modules $V(\omega_\alpha)$ and
$V(\omega_\beta)$ are the simple $\operatorname{SL}_{m+1}$-modules
$\Cbb^{m+1}$ and its dual $(\Cbb^{m+1})^*$, respectively. Let denote
by $Q$ the quadratic form on $\Cbb^{m+1}\oplus(\Cbb^{m+1})^*$
defined by $Q(u,u^*)=\langle u^*,u\rangle$. Then $Q$ is invariant
under the action of $\operatorname{SL_{m+1}}$. Moreover
$Q(v_{\omega_\alpha}+v_{\omega_\beta})=0$, so that $X^1$ is a
subvariety of the quadric ($Q=0$) in
$\Pbb(\Cbb^{m+1}\oplus(\Cbb^{m+1})^*)=\Pbb(\Cbb^{2m+2})$.

We complete the proof by computing the dimension of $X^1$: \begin{equation}\label{dim}
\operatorname{dim}\,X^1=\operatorname{dim}\,G/H=1+\operatorname{dim}\,G/P=1+\sharp(R^+\backslash R_I^+)
\end{equation}
where $R^+$ is the set of positive roots of $(G,B)$ and $R_I^+$ is the set of positive roots generated by simple roots of $I$.
So $\operatorname{dim}\,X^1=\operatorname{dim}\,Q^{2m}=2m$ and $X^1=Q^{2m}$.
\end{proof}

\begin{prop}\label{cas2}
In Case 2, $X^1$ is isomorphic to the grassmannian
$\operatorname{Gr}(i+1,m+2)$.
\end{prop}

\begin{proof}
The fundamental $\operatorname{SL_{m+1}}$-modules are exactly the

$$\begin{array}{c} V(\omega_i)=\bigwedge^i\Cbb^{m+1}
\end{array}$$
 and a highest weight vector of $V(\omega_i)$ is $e_1\wedge\cdots\wedge e_i$ where $e_1,\ldots,e_{m+1}$ is a basis of $\Cbb^{m+1}$.

 We have
$$\xymatrix{
    X\, \ar@{^{(}->}[r]  & \,\Pbb(\bigwedge^i\Cbb^{m+1}\oplus\bigwedge^{i+1}\Cbb^{m+1})\\
    G/H\, \ar@{=}[r] \ar@{^{(}->}[u]  & G.[e_1\wedge\cdots\wedge e_i+e_1\wedge\cdots\wedge e_{i+1}] \ar@{^{(}->}[u]
  }$$
Complete $(e_1,\ldots,e_{m+1})$ to obtain a basis $(e_0,\ldots,e_{m+1})$ of $\Cbb^{m+2}$, then the morphism
$$\begin{array}{ccc}
\bigwedge^i\Cbb^{m+1}\oplus\bigwedge^{i+1}\Cbb^{m+1} & \lra & \bigwedge^{i+1}\Cbb^{m+2}\\
x+y & \lmt & x\wedge e_0 +y
\end{array}$$
 is an isomorphism.
Then $X^1$ is a subvariety of the grassmannian $$\begin{array}{c}
\operatorname{Gr}(i+1,m+2)\simeq
\operatorname{SL}_{m+2}.[e_1\wedge\cdots\wedge
e_i\wedge(e_0+e_{i+1})]\subset\Pbb(\bigwedge^{i+1}\Cbb^{m+1}).\end{array}$$
We conclude by proving that they have the same dimension using
Formula \ref{dim}.
\end{proof}

\begin{prop}
In Case 6, $X^1$ is isomorphic to the spinor variety
$\operatorname{Spin}(2m+1)/P(\omega_m)$.
\end{prop}

\begin{proof}
The direct sum $V(\omega_\alpha)\oplus V(\omega_\beta)$ of the two
half-spin  $\operatorname{Spin}(2m)$-modules is isomorphic to the
spin $\operatorname{Spin}(2m+1)$-module.  Moreover
$v_{\omega_\alpha}+v_{\omega_\beta}$ is in the orbit of a highest
weight vector of the spin $\operatorname{Spin}(2m+1)$-module. Thus
we deduce that $X^1$ is a subvariety of
$\operatorname{Spin}(2m+1)/P(\omega_m)$. We conclude by proving that
they have the same dimension using Formula \ref{dim}.
\end{proof}

\subsection{Non-homogeneous varieties}\label{4}

With the notation of Section \ref{2.2} we prove in this section the following result.

\begin{teo}\label{mainteo}
In Cases 3, 4, 5, 7 and 8 (and only in these cases), $X^1$ is not homogeneous.

Moreover the automorphism group of $X^1$ is
$(\operatorname{SO}(2m+1)\times\Cbb^*)\ltimes V(\omega_m)$,
$(\operatorname{SO}(7)\times\Cbb^*)\ltimes V(\omega_3)$,
$((\operatorname{Sp}(2m)\times\Cbb^*)/\{\pm 1\})\ltimes
V(\omega_1)$, $(\operatorname{F}_4\times\Cbb^*)\ltimes V(\omega_4)$
and $(\operatorname{G}_2\times\Cbb^*)\ltimes V(\omega_1)$
respectively.

Finally, $X^1$ has two orbits under its automorphism group.
\end{teo}

One can remark that in Case 5, Theorem \ref{mainteo} follows from
results of I.~A.~Mihai \cite[Chap.3 and Prop.5.1]{Mi05} combined
with the following result.

\begin{prop}\label{oddgrass}
In Case 5, $X^1$ is isomorphic to the odd symplectic grassmannian
$\operatorname{Gr}_\omega(i+1,2m+1)$.
\end{prop}

\begin{proof}
As in the proof of Proposition \ref{cas2},  $X^1$ is a subvariety of
the odd symplectic grassmannian
$$\operatorname{Gr}_\omega(i+1,2m+1)\simeq\overline{\operatorname{Sp}_{2m+1}.[e_1\wedge\cdots\wedge
e_i\wedge(e_0+e_{i+1})]}\subset\Pbb(\Lambda^{i+1}\Cbb^{m+1}).$$
Again we conclude by proving that they have the same dimension using
Formula \ref{dim} and \cite[Prop 4.1]{Mi05}.
\end{proof}

Now, let $X$ be one of the varieties $X^1$ in Cases 3, 4, 5, 7 and
8.

Then $X$ has three orbits under the action of $G$ (the open orbit
$X_0$ isomorphic to $G/H$ and two closed orbits).  Recall that $X$
can be seen as a subvariety of $\Pbb(V(\omega_\alpha)\oplus
V(\omega_\beta))$.  Let $P_Y:=P(\omega_\alpha)$,
$P_Z:=P(\omega_\beta)$ and denote by $Y$ and $Z$ the closed orbits,
isomorphic to $G/P_Y$ and $G/P_Z$ respectively. (In Case 8 where $G$
is of type $G_2$, we have $\alpha=\alpha_2$ and $\beta=\alpha_1$.)

Let $X_Y$ be  the simple $G/H$-embedding of $X$ with closed orbit
$Y$, we have $X_Y=X_0\cup Y$.  Then, by \cite[Chap.2]{Pa06}, $X_Y$
is a homogeneous vector bundle over $G/P_Y$  in the sense of the
following.
\begin{defi}\label{fibrehom}
Let $P$ be a parabolic subgroup of $G$ and $V$ a $P$-module.  Then
the homogeneous vector bundle $G\times^PV$ over $G/P$ is the
quotient of the product $G\times V$ by the equivalence relation
$\sim$ defined by $$\forall g\in G,\,\forall p\in P, \forall v\in
V,\quad (g,v)\sim(gp^{-1},p.v).$$
\end{defi}
Specifically, $X_Y=G\times^{P_Y}V_Y$ where $V_Y$ is a simple
$P_Y$-module of highest weight $\omega_\beta-\omega_\alpha$, and
similarly, $X_Z=G\times^{P_Z}V_Z$ where $V_Z$ is a simple
$P_Z$-module of highest weight $\omega_\alpha-\omega_\beta$.

Denote by $\Aut(X)$ the automorphism group of $X$ and  $\Aut ^0 (X)$
the connected component of $\Aut (X)$ containing the identity.
\begin{rem}
Observe that $\Aut(X)$ is a linear algebraic group. Indeed $\Aut(X)$
acts on the Picard group of $X$ which equals $\Zbb$ (the Picard
group of a projective spherical variety is free \cite{Br89}). This
action is necessarily trivial. Then $\Aut(X)$ acts on the
projectivization of the space of global sections of a very ample
bundle. This gives a faithful projective representation of
$\Aut(X)$.
\end{rem}

We now complete the proof of Theorem \ref{mainteo} by proving several lemmas.

\begin{lem}\label{fixedorbit}
The closed orbit $Z$ of $X$ is stable under $\Aut ^0 (X)$.
\end{lem}

\begin{proof}
We are going to prove that the normal sheaf $N_Z$ of $Z$ in $X$ has
no nonzero global section. This will imply that $\Aut^0(X)$
stabilizes $Z$, because the Lie algebra $\Lie(\Aut^0(X))$ is the
space of global sections $H^0(X,T_{X})$ of the tangent sheaf $T_X$
 \cite[Chap.2.3]{Ak95} and we have the following exact
sequence
$$0\lra T_{X,Z}\lra T_{X}\lra N_Z\lra 0$$
where $T_{X,Z}$ is the subsheaf of $T_X$ consisting of vector fields
that vanish along $Z$.  Moreover $H^0(X,T_{X,Z})$ is the Lie algebra
of the subgroup of $\Aut^0(X)$ that stabilizes $Z$.

The total space of $N_Z$ is the vector bundle $X_Z$. So using the
Borel-Weil theorem \cite[4.3]{Ak95},  $H^0(G/P_Z,N_Z)=0$ if and only
if the smallest weight of $V_Z$ is not antidominant. The smallest
weight of $V_Z$ is  $w_0^\beta(\omega_\alpha-\omega_\beta)$ where
$w_0^\beta$ is the longest element of $W_{S\backslash\{\beta\}}$.
Let $\gamma\in S$, then $$\langle
w_0^\beta(\omega_\alpha-\omega_\beta),\check\gamma\rangle=\langle\omega_\alpha-\omega_\beta,w_0^\beta(\check\gamma)\rangle.$$
If $\gamma$ is different from $\beta$  then
$w_0^\beta(\check\gamma)=-\check\delta$ for some $\delta\in
S\backslash\{\beta\}$. So we only have to compute
$w_0^\beta(\check\beta)$.

In Case 3, $\beta=\alpha_m$, so $w_0^\beta$ maps $\alpha_i$ to $-\alpha_{m-i}$ for all $i\in\{1,\ldots,m-1\}$.
 Here, $\omega_\beta=\frac{1}{2}(\alpha_1+2\alpha_2+\cdots+m\alpha_m)$.
 Then, using the fact that $w_0^\beta(\omega_\beta)=\omega_\beta$, we have  $w_0^\beta(\beta)=\alpha_1+\cdots+\alpha_m$
 so that $w_0^\beta(\check\beta)=2(\check\alpha_1+\cdots+\check\alpha_{m-1})+\check\alpha_m$ and $\langle\omega_\alpha-\omega_\beta,w_0^\beta(\check\gamma)\rangle=1$
 (because $\alpha=\alpha_{m-1}$).

 The computation of $w_0^\beta(\check\beta)$ in the other cases is similar and left to the reader.
 In all four cases,$\langle\omega_\alpha-\omega_\beta,w_0^\beta(\check\beta)\rangle>0$ (this number equals 1 in Cases 3, 4, 5, 7 and 2 in Case 8).
This proves that $w_0^\beta(\omega_\alpha-\omega_\beta)$ is not antidominant.
\end{proof}

\begin{rem}\label{sectionNY}
Using the Borel-Weil theorem we can also compute $H^0(G/P_Y,N_Y)$ by
the same method.  We find that in Cases 3, 4, 5, 7 and 8, this
$G$-module is isomorphic to the simple $G$-module $V(\omega_m)$,
$V(\omega_3)$, $V(\omega_1)$, $V(\omega_4)$ and $V(\omega_1)$
respectively.
\end{rem}

We now prove the following lemma.

\begin{lem}\label{lemaut}
In Cases 3, 4, 5, 7 and 8, $\Aut^0(X)$ is
$(\operatorname{SO}(2m+1)\times\Cbb^*)\ltimes V(\omega_m)$,
$(\operatorname{SO}(7)\times\Cbb^*)\ltimes V(\omega_3)$,
$(\operatorname{Sp}(2m)\times\Cbb^*)\ltimes V(\omega_1)$,
$(\operatorname{F}_4\times\Cbb^*)\ltimes V(\omega_4)$ and
$(\operatorname{G}_2\times\Cbb^*)\ltimes V(\omega_1)$ respectively.
\end{lem}
\begin{rem}
By Remark \ref{remP/H}, we already know that the action of $G$ on
$X$ extends to an action of $G\times\Cbb^*$($\simeq G\times P/H$).
Moreover $\tilde{C}:=\{(c,c^{-1}H)\in C\times P/H\}\subset
G\times\Cbb^*$ acts trivially on $X$.
\end{rem}
\begin{proof}

Let $\pi:\tilde{X}\lra X$ be the blowing-up of $Z$ in $X$. Since $Z$
and $X$ are smooth, $\tilde{X}$ is smooth; it is also a projective
$G/H$-embedding. In fact $\tilde{X}$ is the projective bundle
$$\phi:\,G\times^{P_Y}\Pbb(V_Y\oplus\Cbb)\lra G/P_Y$$ where $P_Y$
acts trivially on $\Cbb$ (it is the unique projective
$G/H$-embeddings with unique color $\alpha$). Moreover the
exceptional divisor $\tilde{Z}$ of $\tilde{X}$ is $G/P$.

Let us remark that $\Aut^0(\tilde{X})$ is isomorphic to $\Aut^0(X)$.
Indeed, it contains $\Aut^0(X)$ because $Z$ is stable under the
action of $\Aut^0(X)$. Conversely, we know, by a result of
A.~Blanchard, that $\Aut^0(\tilde{X})$ acts on $X$ such that $\pi$
is equivariant \cite[Chap.2.4]{Ak95}.

Now we are going to compute $\Aut^0(\tilde{X})$.
Observe that $H^0(G/P_Y,N_Y)$ acts on $\tilde{X}$ by translations on the fibers of $\phi$:
$$\forall s\in H^0(G/P_Y,N_Y),\, \forall (g_0,[v_0,\xi])\in G\times^{P_Y}\Pbb(V_Y\oplus\Cbb),\,
s.(g_0,[v_0,\xi])=(g_0,[v_0+\xi v(g_0),\xi])$$ where $v(g_0)$ is the element of $V_Y$ such that $s(g_0)=(g_0,v(g_0))$.

Then the group $((G\times\Cbb^*)/\tilde{C})\ltimes H^0(G/P_Y,N_Y)$ acts effectively on $\tilde{X}$
(the semi-product is defined by $((g',c'),s').((g,c),s)=((g'g,c'c),c'g's+s')$).
In fact we are going to prove that $$\Aut^0(\tilde{X})=((G\times\Cbb^*)/\tilde{C})\ltimes H^0(G/P_Y,N_Y).$$

By \cite[Chap.2.4]{Ak95} again, we know that $\Aut^0(\tilde{X})$
exchanges the fibers of $\phi$ and induces an automorphism of
$G/P_Y$. Moreover we have $\Aut^0(G/P_Y)=G/C$  in our four cases
\cite[Chap.3.3]{Ak95}. So we have an exact sequence
$$0\lra A\lra \Aut^0(\tilde{X})\lra G/C\lra 0$$
where $A$ is the set of automorphisms which stabilize each  fiber of
the projective bundle $\tilde{X}$. In fact $A$ consists of affine
transformations in fibers.

Then $H^0(G/P_Y,N_Y)$ is the subgroup of $A$ consisting of
translations. Let $A_0$ be the subgroup of $A$ consisting of linear
transformations in fibers. Then $A_0$ fixes $Y$ so that $A_0$ acts
on the blowing-up $\tilde{\tilde{X}}$ of $Y$ in $\tilde{X}$.
Moreover $\tilde{\tilde{X}}$ is a $\Pbb^1$-bundle over $G/P$, it is
in fact the unique $G/H$-embedding without colors \cite[Ex.1.13
(2)]{Pa06}. As before we know that $\Aut^0(\tilde{\tilde{X}})$
exchanges the fibers of that $\Pbb^1$-bundle and induces an
automorphism of $G/P$.
Moreover we have $\Aut^0(G/P)=G/C$ and then $\Aut^0(\tilde{\tilde{X}})=(G\times\Cbb^*)/\tilde{C}$. We deduce that $A_0=\Cbb^*$.\\

We complete the proof of Lemma \ref{lemaut} by Remark
\ref{sectionNY}.

\end{proof}

To complete the proof of Theorem \ref{mainteo} we prove the following lemma.

\begin{lem}
The automorphism group of $X$ is connected.
\end{lem}

\begin{proof}
Let $\phi$ be an automorphism of $X$. We want to prove that $\phi$
is in $\Aut^0(X)$.   But $\phi$ acts by conjugation on $\Aut^0(X)$.
Let $L$ be a Levi subgroup of $\Aut^0(X)$. Then $\phi^{-1}L\phi$ is
again a Levi subgroup of $\Aut^0(X)$. But all Levi subgroups are
conjugated in $\Aut^0(X)$. So we can suppose, without loss of
generality, that $\phi$ stabilizes $L$.

Then $\phi$ induces an automorphism of the direct product of
$\Cbb^*$ with a simple group $G$ of type $B_m$, $C_m$, $F_4$ or
$G_2$ (Lemma \ref{lemaut}).  It also induces an automorphism of $G$
which is necessarily an inner automorphism (because there is no
non-trivial automorphism of the Dynkin diagram of $G$). So we can
assume now that $\phi$ commutes with all elements of $G$.

Then $\phi$ stabilizes the open orbit $G/H$ of $X$. Let $x_0:=H\in
G/H\subset X$ and $x_1:=\phi(x_0)\in G/H$. Since $\phi$ commutes
with the elements of $G$, the stabilizer of $x_1$ is also $H$. So
$\phi$ acts on $G/H$ as an element of $N_G(H)/H=P/H\simeq \Cbb^*$
(where $N_G(H)$ is the normalizer of $H$ in $G$).  Then $\phi$ is an
element of $\Cbb^*\subset\Aut^0(X)$.
\end{proof}

\begin{rem}
In Case 5,  we recover the result of I.~Mihai:
$\Aut(X)=((\operatorname{Sp}(2m)\times\Cbb^*)/\{\pm
1\})\ltimes\Cbb^{2m}$.
\end{rem}

\subsection{First step in the proof of Theorem \ref{caracterisation2}}\label{5}

In this section, we prove the following result.

\begin{teo}\label{caracterisation}
Let $X$ be a smooth projective variety with Picard number~1 and put
$\mathbf{G}:=\Aut^0(X)$. Suppose that $\mathbf{G}$ is not
semi-simple and that $X$ has two orbits under the action of
$\mathbf{G}$. Denote by $Z$ the closed orbit.

Then the codimension of $Z$ is at least $2$.

Suppose in addition that the blowing-up of $X$ along $Z$ has also
two orbits under the action of $\mathbf{G}$.  Then $X$ is one of the
varieties $X^1$ obtained in Cases 3, 4, 5, 7 and 8 of Theorem
\ref{classif}.
\end{teo}

\begin{rem}
The converse implication holds by Theorem \ref{mainteo}.
\end{rem}

We will use a result of D.~Akhiezer on two-orbits varieties with
codimension one closed orbit.

\begin{lem}[Th.1 of \cite{Ah83}]\label{Akhiezer}

Let $X$ a smooth complete variety with an effective action of the
connected linear non semi-simple group $\mathbf{G}$. Suppose that
$X$ has two orbits under the action of $\mathbf{G}$ and that the
closed orbit is of codimension $1$.
Let $G$ be a maximal semi-simple subgroup of $\mathbf{G}$.\\

Then there exist a parabolic subgroup $P$ of $G$ and a $P$-module $V$ such that:\\
(i) the action of $P$ on $\Pbb(V)$ is transitive;\\
((ii) there exists an irreducible $G$-module $W$ and a surjective $P$-equivariant morphism $W\lra V$;)\\
(iii) $X=G\times^P\Pbb(V\oplus\Cbb)$ as a $G$-variety.\\
\end{lem}

\begin{rem}
It follows from (i) that $V$ is an (irreducible) horospherical
$L$-module of rank $1$, where $L$ is a Levi subgroup of $G$. This
implies that $X$ is a horospherical $G$-variety of rank $1$.

(The $G$-module $W$ is the set of global sections of the vector
bundle $G\times^P V\lra G/P$.)
\end{rem}

\begin{proof}[Proof of Theorem \ref{caracterisation}]
If $Z$ is of codimension $1$,  Lemma \ref{Akhiezer} tells us that
$X$ is a horospherical $G$-variety. Moreover the existence of the
irreducible $G$-stable divisor $Z$ tells us that one of the two rays
of the colored fan of $X$ has no color \cite[Chap 1]{Pa06}. Then $X$
satisfies the condition (i) of Theorem \ref{thesproj} and
$X=\Pbb(V\oplus\Cbb)$. Since $X$ is not homogeneous we conclude that
the codimension of $Z$ is at least~2.

Denote by $\tilde{X}$ the blowing-up of $X$ along $Z$. Then
$\tilde{X}$ is a horospherical $G$-variety by Lemma \ref{Akhiezer}.
Moreover $X$ and $\tilde{X}$ have the same open $G$-orbit, so that
$X$ is also a horospherical $G$-variety. We conclude by the
description of smooth projective horospherical varieties with Picard
number~1 obtained in the preceding sections .
\end{proof}

\section{On some two-orbits varieties under the action of a semi-simple group}\label{section2}

The aim of this section is to prove Theorem \ref{caracterisation2}.

We keep the notation of the first paragraph of Section \ref{1}. In
all this section, $G$ is a semi-simple group and all varieties are
spherical of rank one (but not horospherical).

We will often use the following result, that can be deduced from the
local structure of spherical varieties \cite[Chap.1]{Br89}.

\begin{prop}\label{structlocale}
Let $X$ be a (spherical) tow-orbits $G$-variety, with closed orbit
isomorphic to $G/Q$ (choose $Q\supset B\supset T$ such that $X$ has
an open $B$-orbit). Let $L(Q)$ be the Levi subgroup of $Q$
containing $T$ and $Q^-$ the parabolic subgroup of $G$ opposite to
$Q$. Then there exists a affine $L(Q)$-subvariety $V$ of $X$ with a
fixed point and an affine open subvariety $X_0$ of $X$, such that
the natural morphism from $R_u(Q^-)\times V$ to $X$ is an
isomorphism into $X_0$. Moreover $V$ is spherical of rank one.
\end{prop}

\subsection{First results}\label{genres}

Let us first prove the first part of Theorem \ref{caracterisation2}.

\begin{prop}\label{codim>2}
Let $X$ be a non-homogeneous projective two-orbits variety with
Picard number~1. Then the codimension of the closed orbit $Z$ is at
least $2$.
\end{prop}

\begin{proof}
Suppose that $Z$ is of codimension~1. This implies that $X$ is
smooth.

Let $G/H$ be a homogenous space isomorphic to the open orbit of $X$
and $P$ be a parabolic subgroup of $G$ containing $H$ and minimal
for this property.

Then $P$ satisfies one of the two following conditions \cite[proof
of Th.5]{Ah83}:
\begin{itemize}
\item $R(P)\subset H\subset P$,
\item $H$ contains a Levi subgroup of $P$.
\end{itemize}

Moreover, if $P=G$, then $X$ is unique and homogeneous under its full automorphism group \cite[Th.4]{Ah83}. So $P$ is a proper subgroup of $G$.\\

We may assume, without loss of generality, that $G$ is simply
connected.

We have the following exact sequence.
$$\Zbb[Z] \lra \operatorname{Pic}(X) \lra \operatorname{Pic}(G/H)
\lra 0$$ so that the Picard group  $\operatorname{Pic}(G/H)$ of
$G/H$ is finite. But, by \cite[Th.2.2]{Ti06},
$\operatorname{Pic}(G/H)$ is isomorphic to the group
$\mathfrak{X}(H)$ of characters of $H$ . This contradicts the fact
that $H$ contains the radical or a Levi subgroup of a proper
parabolic subgroup of $G$. Indeed we have one of the two following
commutative diagrams (respectively when $R(P)\subset H$ and
$L(P)\subset H$),
$$\xymatrix{
\mathfrak{X}(P)\, \ar[r]\ar[rd] & \mathfrak{X}(H) \ar[d] &  &\mathfrak{X}(P)\, \ar[r]\ar[rd] & \mathfrak{X}(H) \ar[d] & \\
& \mathfrak{X}(R(P)) & & & \mathfrak{X}(L(P)) }$$ where the maps are
restrictions and $L(P)$ is a Levi subgroup of $P$. Moreover the maps
$\mathfrak{X}(P)\lra\mathfrak{X}(R(P))$ and
$\mathfrak{X}(P)\lra\mathfrak{X}(L(P))$ are injective. Then
$\mathfrak{X}(P)=0$ so that $P=G$.
\end{proof}

From now on, let $X$ be a two-orbits variety satisfying (*)
(Definition \ref{etoile}). Let $G/H$ be a homogenous space
isomorphic to the open orbit of $X$ and $P$ be a parabolic subgroup
of $G$ containing $H$ and minimal for this property. Then, we still
have the two cases of Proposition \ref{codim>2}: $R(P)\subset
H\subset P$ or $H$ contains a Levi subgroup of $P$.

Denote by $\pi:\tilde{X}\lra X$ the blow-up of $X$ along $Z$. The
variety $\tilde{X}$ is the unique $G/H$-embedding ({\it i.e.} a
normal $G$-variety with an open orbit isomorphic to $G/H$) with
closed orbit of codimension~1.

Note that $\tilde{X}$ cannot be homogeneous under its full
automorphism group, because of a
result of A.~Blanchard \cite[Chap.2.4]{Ak95}.\\

Let us remark also that $X$ is the unique $G/H$-embedding satisfying
(*). Indeed, a spherical homogeneous space of rank one cannot have
two different projective $G/H$-embeddings with Picard number one. In
fact, if it exists, the projective embedding with Picard number one
is the $G/H$-embeddings associated to the unique complete colored
fan with all possible colors (see, for example \cite{Kn91} for the
classification of $G/H$-embeddings, and \cite{Br89} for the
description of the Picard group of spherical varieties).

To conclude this section, let us prove the following.
\begin{lem}
$P$ is a maximal parabolic subgroup of $G$.
\end{lem}

\begin{proof}
Let us assume, without loss of generality, that $G$ is simply
connected. Since $Z$ has codimension at least~2, the Picard groups
of $G/H$ and $X$ are the same. The argument of the proof of Prop
\ref{codim>2} implies that $\mathfrak{X}(P)=\Zbb$ so that $P$ is
maximal.
\end{proof}

\subsection{When $R(P)\subset H$}\label{RPH}

We suppose in all this section that $R(P)\subset H$. Let us remark
that $P/H$ is isomorphic to $(P/R(P))/(H/R(P))$ and is still
spherical of rank one. Then we deduce the following lemma from
Theorem 4 of \cite{Ah83}.

\begin{lem}\label{3cas}
Suppose that there exist a projective $G/H$-embedding with  Picard
number~1 and that $R(P)\subset H\subset P\neq G$. Then $P/H$ is
isomorphic to one of the following homogeneous spaces $G'/H'$ where
$G'$ is a quotient of a normal subgroup of $P/R(P)$ by a finite
central subgroup:
\begin{center}
\begin{tabular}{|c||c|c|c|c|}
\hline
 & $G'$ & $H'$ & $X'$ & $Q'$\\
\hline
\hline
1a & $\operatorname{SO(n)},\,n\geq 4$ & $\operatorname{SO(n-1)}$ & $Q^{n-1}$ & $P(\omega_{\alpha_1})$ or $B$ if $n=4$\\
\hline
1b & $\operatorname{SO(n)}/C',\,n\geq 4$ & $\operatorname{S(O(1)}\times\operatorname{O(n-1))}/C'$ & $\Pbb^{n-1}$ & $P(\omega_{\alpha_1})$ or $B$ if $n=4$\\
\hline
2 & $\operatorname{Sp(2n)}/C',\,n\geq 2$ & $(\operatorname{Sp(2n-2)}\times \operatorname{Sp(2)})/C'$ & $\operatorname{Gr}(2,2n)$ & $P(\omega_{\alpha_2})$\\
\hline
3a & $\operatorname{Spin(7)}$ & $\operatorname{G_2}$ & $Q^7$ & $P(\omega_{\alpha_3})$\\
\hline
3b & $\operatorname{SO(7)}$ & $\operatorname{G_2}$ & $\Pbb^7$ & $P(\omega_{\alpha_3})$\\
\hline
\end{tabular}
\end{center}

~\\ Here we denote by $C'$ the center of $G'$. The variety $X'$
denotes the unique projective $G'/H'$-embedding.  Its closed orbit
is of codimension 1 and isomorphic to $G'/Q'$ where $Q'$ is the
parabolic subgroup given in the last column. We denote by $Q^{n-1}$
the quadric of dimension $n-1$ in $\Pbb^n$ and by
$\operatorname{Gr}(2,2n)$ the grassmannian of planes in $\Cbb^{2n}$.

Let us remark that, in case 1a and 1b when $n=4$, the group $G'$ is
of type $A_1\times A_1$. This is the only case where $G'$ is not
simple.

Moreover, the closed orbit $Z$ of $X$ is isomorphic to $G/Q$, where
$Q$ is a parabolic subgroup of $G$ containing $B$ such that
$P/(P\cap Q)\simeq G'/Q'$.
\end{lem}

\begin{proof}
We apply \cite[Th.4]{Ah83} to the pair $(P/R(P), H/R(P))$ (note that
$H/R(P)$ is reductive because $P$ is a minimal parabolic subgroup
containing $H$).  So $P/H$ is isomorphic to a homogeneous space
listed in \cite[Table 2]{Ah83}. Since $P$ is a proper parabolic
subgroup of $G$, the group $G'$ cannot be $\operatorname{G_2}$ or
$\operatorname{F_4}$.

Moreover we have already seen that the restriction morphism
$\mathfrak{X}(H)\lra\mathfrak{X}(R(P))$ is injective, hence
$\mathfrak{X}(H')$ is finite. So the cases
$(G'=\operatorname{PSL(n+1)},\,H'=\operatorname{GL(n)})$,
$(G'=\operatorname{SO(3)},\, H'=\operatorname{SO(2)})$ and
$(G'=\operatorname{SO(3)},\,
H'=\operatorname{S(O(1)}\times\operatorname{O(2)}))$ cannot happen.

For the last statement, remark that $\tilde{X}=G\times^PX'$ (see
Definition \ref{fibrehom}), the closed orbit of $\tilde{X}$ is
$G/(P\cap Q)$ and the closed orbit of $X'$ is $G'/Q'\simeq P/(P\cap
Q)$. Note that the blow-up $\pi:\tilde{X}\lra X$ is proper is
$G$-equivariant and it sends the closed orbit of $\tilde{X}$ to the
closed orbit $Z$ of $X$, so that $Z$ is isomorphic to $G/Q$, $G/P$
or to a point. By Proposition \ref{structlocale}, a two-obits
variety whose closed orbit is a point must be affine, so that $Z$
cannot be point. Now if $Z$ is isomorphic to $G/P$, again by
Proposition \ref{structlocale}, there exists an open subvariety of
$X$ isomorphic to the product of an open subvariety of $G/P$ and a
closed subvariety $S$ of $X$ spherical under the action of a Levi
subgroup of $P$. Moreover $S$ is necessarily isomorphic to a
projective $G'/H'$-embedding with closed orbit a point, that gives
us a contradiction.
\end{proof}

Let us now translate the fact that $X$ is smooth.

\begin{lem}\label{smooth}
Suppose that there exists a smooth $G/H$-embedding $X$ with closed
orbit $G/Q$.

Denote by $Q^-$ the parabolic subgroup of $G$ opposite to $Q$ and by
$L(Q):=Q\cap Q^-$ the Levi subgroup of $Q$ containing $T$.

Then, there is an affine open subvariety of $X$ that is isomorphic
to the product of $R_u(Q^-)$ and a $L(Q)$-module $V$ included in
$X$.

Moreover, $V$ is either a fundamental module of a group of type
$A_n$ and of highest weight $\omega_1$  or $\omega_n$, or a
fundamental module of a group of type $C_n$ and of highest weight
$\omega_1$ (modulo the action of the center of $Q$).
\end{lem}

\begin{rem}\label{structurelocaleX'}
For the variety $X'$, one can directly compute that there is an
affine open subvariety of $X'$ that is isomorphic to the product of
$R_u(Q'^-)$ and a line where $L(Q')$ acts respectively with weight:

$-\omega_1$ in Case 1a ($-\omega_1-\omega_1$ when $n=4$);

$-2\omega_1$ in Case 1b ($-2\omega_1-2\omega_1$ when $n=4$);

$-\omega_2$ in Case 2;

$-\omega_3$ in Case 3a;

$-2\omega_3$ in Case 3b.

This remark will be used to prove that, in some cases, $X$ cannot be smooth.
\end{rem}

\begin{proof}
By Proposition \ref{structlocale}, there exists a smooth affine
$L(Q)$-subvariety $V$ of $X$ with a fixed point and an affine open
subvariety of $X$ that is isomorphic to $R_u(Q^-)\times V$. Moreover
$V$ is spherical of rank one.

But if $L$ is a reductive connected algebraic group acting on a
smooth affine variety $V$ with a fixed point, then $V$ is a
$L(Q)$-module \cite[Cor.6.7]{PV94}. So $V$ is a $L(Q)$-module.

Moreover $V$ has two $L(Q)$-orbits, so that it is horospherical.
Then, the last part of the lemma follows from \cite[Lem.2.13]{Pa06}.
\end{proof}

\begin{cor}\label{corsmooth}
Suppose that there exists a $G/H$-embedding $X$ satisfying (*).
Suppose also that $Q$ is maximal (this holds in all cases except
Case 1 of Lemma \ref{3cas} with $n=4$). Let $i$ and $j$ be the
indices such that $P=P(\omega_i)$ and $Q=P(\omega_j)$.

Let us define $\Gamma_i^j$ (respectively $\Gamma_j^i$) to be the
connected component of the subgraph
$\Gamma_{S\backslash\{\alpha_i\}}$ (respectively
$\Gamma_{S\backslash\{\alpha_j\}}$) of the Dynkin diagram of $G$
containing $\alpha_j$ (respectively $\alpha_i$).

Then $\Gamma_i^j$ is of one of the following types:\\
$B_n$, $n\geq 2$, with first vertex $\alpha_j$, \\
$D_n$, $n\geq 3$, with first vertex  $\alpha_j$,
\\
$C_n$, $n\geq 2$, with second vertex $\alpha_j$,
\\
$B_3$, with third vertex $\alpha_j$.

And $\Gamma_j^i$ is of type $A_n$ with first (or last) vertex
$\alpha_i$ or of type $C_n$ with first vertex $\alpha_i$.
\end{cor}

\begin{proof}
The first part of the corollary follows from Lemma \ref{3cas}.

We deduce the second part from Lemma \ref{smooth}. Recall that
$G/H=G\times^PP/H$ and remark that Lemma \ref{smooth} applied to the
$P/H$-embedding $X'$ says that there is an open affine subset of
$X'$ that is isomorphic to the product of $R_u((P\cap Q)^-)$ and a
$L((P \cap Q)^-)$-module of dimension~1. Then the highest weight of
the $L(Q)$-module $V$ must be a character of $P\cap Q$. We conclude
by the last statement of Lemma \ref{smooth} applied to $X$.
\end{proof}

When $G$ is simple, applying Corollary \ref{corsmooth}, we are able
to give a first list of possible homogeneous spaces admitting an
embedding that satisfies (*).

\begin{lem}\label{touslescas}
Suppose that there exists a $G/H$-embedding $X$ satisfying (*).
Suppose that $G$ is simple and $R(P)\subset H$, then $(G,P,Q)$ is
one of the following:\\

(a) $(A_4,P(\omega_1),P(\omega_3))$ or, that is the same,
$(A_4,P(\omega_4),P(\omega_2))$;\\

(b) $(B_n,P(\omega_i),P(\omega_{i+1}))$, with $n\geq 3$ and $1\leq
i\leq n-2$,

 ~or $(D_n,P(\omega_i),P(\omega_{i+1}))$, with $n\geq 4$ and $1\leq
i\leq n-3$,

 ~or $(D_n,P(\omega_{n-2}),P(\omega_{n-1})\cap
P(\omega_n))$, with $n\geq 3$;\\

(c) $(B_4,P(\omega_4),P(\omega_2))$;\\

(d) $(B_4,P(\omega_1),P(\omega_4))$;\\

(e) $(B_3,P(\omega_2),P(\omega_1)\cap P(\omega_3))$;\\

(f) $(C_n,P(\omega_1),P(\omega_3))$ with $n\geq 3$;\\

(g) $(C_3,P(\omega_2),P(\omega_1)\cap P(\omega_3))$;\\

(h) $(F_4,P(\omega_1),P(\omega_3))$;\\

(i) $(F_4,P(\omega_4),P(\omega_1))$;\\

(j) $(F_4,P(\omega_4),P(\omega_3))$.\\
\end{lem}






In the next Lemma, we consider the case where $G$ is not simple.

\begin{lem}\label{casnotsimple}
Suppose that there exists a $G/H$-embedding $X$ satisfying (*).
Suppose that $G$ is not simple, acts faithfully on $G/H$ and
$R(P)\subset H$. Then $(G,P,Q)$ is one of the following:\\

(a') $(A_n\times A_1,P(\omega_2)\times A_1,P(\omega_1)\times
P(\omega_1))$ or, that is the same, $(A_n\times
A_1,P(\omega_{n-1})\times A_1,P(\omega_n)\times P(\omega_1))$, with $n\geq 2$;\\

(b') $(B_n\times A_1,P(\omega_{n-1})\times A_1,P(\omega_n)\times
P(\omega_1))$ with $n\geq 3$;\\

(c') $(C_n\times A_1,P(\omega_{n-1})\times A_1,P(\omega_n)\times
P(\omega_1))$ with $n\geq 2$; \\

(d') $(C_n\times A_1,P(\omega_2)\times A_1,P(\omega_1)\times
P(\omega_1))$ with $n\geq 2$;\\

(e') $(G_2\times A_1,P(\omega_1)\times A_1,P(\omega_2)\times
P(\omega_1))$;\\

(f') $(G_2\times A_1,P(\omega_2)\times A_1,P(\omega_1)\times
P(\omega_1))$;\\
\end{lem}

\begin{proof}
We may assume that $G=G_{i_1}\times\cdots\times G_{i_k}$ where
$k\geq 2$ and $G_{i_1},\ldots,G_{i_k}$ are simple groups. Since $P$
is maximal, we may assume also that $P=P_{i_1}\times
G_{i_2}\times\ldots G_{i_k}$ and then
$P/R(P)=P_{i_1}/R(P_{i_1})\times G_{i_2}\times\ldots G_{i_k}$.

Moreover $P/H$ is isomorphic to one of the homogeneous space of
Lemma \ref{3cas}. Since $G$ is not simple and acts faithfully on
$G/H$, $H$ cannot contain the subgroup $G_{i_2}\times\cdots\times
G_{i_k}$ so that $P/H$ is isomorphic to
$(\operatorname{SL(2)}\times\operatorname{SL(2)})/\operatorname{SL(2)}$
or
$(\operatorname{PSL(2)}\times\operatorname{PSL(2)})/\operatorname{PSL(2)}$
({\it i.e.} Case 1a and 1b of Lemma \ref{3cas} with $n=4$,
respectively). One can deduce that $k=2$, that $G_{i_2}$ and a
normal simple subgroup of $P_{i_1}/R(P_{i_1})$ are of type $A_1$.
\end{proof}

\begin{rem}
In each case of Lemmas \ref{touslescas} and \ref{casnotsimple} , there exist at most one homogeneous space $G/H$ and
one $G/H$-embedding satisfying (*). Indeed, suppose that there exist
 two varieties $X_a$ and respectively $X_b$ that satisfy the same case of Lemma
\ref{touslescas} and Case 1a and respectively Case 1b (or Case 3a
and respectively Case 3b) of Lemma \ref{3cas}. Then $X_a$ is a
double cover of $X_b$ ramified along the closed orbit $G/Q$ (because
the quadric $Q^n$ is a ramified double cover of $\Pbb^n$). By the
purity of the branch locus, both $X_a$ and $X_b$ cannot be smooth,
that contradicts the hypothesis of the beginning.
\end{rem}

In the next three subsections we study separately all the cases
enumerated in Lemmas \ref{touslescas} and \ref{casnotsimple}.

\subsubsection{Non-homogeneous varieties}
Let us first define two varieties as follows:

\begin{defi}\label{defiX1}
Let $G=F_4$, $P=P(\omega_1)$ and $Q=P(\omega_3)$.

Denote by $L(P)$ the Levi subgroup of $P$ containing $T$ and by $\omega_i'$ the fundamental
weights of $(P/R(P), B/R(P))$ (here $P/R(P)$ is of type $C_3$). Let
$V$ be the $G$-module $V(\omega_1)\oplus V(\omega_3)$, $V'$ be the
sub-$L(P)$-module of $V$ generated by $v_{\omega_3}$ and $\Cbb
.v_{\omega_1}$ be the line of $V$ generated by $v_{\omega_1}$.
Remark that $V'$ is the fundamental $L(P)$-module of weight
$\omega_3=\omega_2'+\omega_1$.

Let $X'$ be the two-orbits $P/R(P)$-variety of Case 2 of Lemma
\ref{3cas}, included in $\Pbb(\Cbb .v_{\omega_1}\oplus V')$ ($X'$ is isomorphic to the grassmannian $\operatorname{Gr}(2,6)$).

One can now define $\mathbf{X_1}:=\overline{G.X'}\subset \Pbb(V)$.
\end{defi}

\begin{defi}\label{defiX2}
Let $G=G_2\times \operatorname{PSL(2)}$ acting on
$\Pbb(\Im(\Obb)\otimes\Cbb^2)$ where $\Im(\Obb)$ is the
non-associative algebra consisting of imaginary octonions and $G_2$
is the group of automorphism of $\Im(\Obb)$. (See proof of
Proposition \ref{casiv} for more details about octonions).

Let $(e_1,e_2)$ be a basis of $\Cbb^2$ and let $z_1,z_2$ be two
elements of $\Im(\Obb)$ such that $z_1^2=z_2^2=z_1z_2=0$. Define
$$x:=[z_1\otimes e_1+z_2\otimes
e_2]\in\Pbb(\Im(\Obb)\otimes\Cbb^2)$$ and $\mathbf{X_2}$ as the
closure of $G.x$ in $\Pbb(\Im(\Obb)\otimes\Cbb^2)$.
\end{defi}

\begin{prop}\label{casf}
The variety $\mathbf{X_1}$ satisfies (*) in Case (h) of
Lemma \ref{touslescas}.
\end{prop}

\begin{proof}
 By construction, the open orbit $G/H$ of $\mathbf{X_1}$
 satisfies $R(P)\subset H \subset P$ such that $P/H$ is as Case
2 of Lemma \ref{3cas}. Moreover, its closed orbit is
$G.[v_{\omega_3}]$ and isomorphic to $G/Q$.

Now let us check that $\mathbf{X_1}$ is smooth. By Remark
\ref{structurelocaleX'}, there exists an affine open subset of $X'$
that is a product of $R_u((P\cap Q)^-)$ and a line where $P\cap Q$
acts with weight $-\omega_2'=\omega_1-\omega_3$. Then there exists
an open subset of $\mathbf{X_1}$ that is a product of $R_u(Q^-)$ and the
irreducible horospherical $L(Q)$-module of highest weight
$\omega_1-\omega_3$.
\end{proof}

\begin{prop}\label{casf'}\label{P5}
The variety $\mathbf{X_2}$ satisfies (*) in Case (f') of Lemma \ref{casnotsimple}.
\end{prop}

\begin{proof}
Let $x$ as in Definition \ref{defiX2} and $H:=\operatorname{Stab}_Gx$.
Then $H$ is included in the maximal parabolic subgroup $P$ of $G$
that stabilizes the plane generated by $z_1$ and $z_2$. Moreover
$R(P)\subset H$ and $P/H$ is isomorphic to
$(\operatorname{PSL(2)}\times\operatorname{PSL(2)})/\operatorname{PSL(2)}$
where $\operatorname{PSL(2)}$ is included in
$\operatorname{PSL(2)}\times\operatorname{PSL(2)}$ as follows:
$$\begin{array}{ccc}
\operatorname{PSL(2)} & \lra & \operatorname{PSL(2)}\times\operatorname{PSL(2)}\\
A & \lmt & (A,^tA^{-1}).
\end{array}$$

Now let us remark that the closed orbit of the closure $X$ of $G.x$
in $\Pbb(\Im(\Obb)\otimes\Cbb^2)$ is $G.[z_1\otimes e_1]$ that is
isomorphic to $G/Q$ where $Q=P(\omega_1)\times P(\omega_0)$ (to
avoid confusion, we denote here by $\omega_0$ the fundamental weight
of $\operatorname{PSL(2)}$).

Let us now check that $X$ is smooth. By Remark
\ref{structurelocaleX'}, there exists an affine open subset of $X'$
that is a product of $R_u((P\cap Q)^-)$ and a line where $P\cap Q$
acts with weight $\omega_2-2\omega_1-2\omega_0$ (where $\omega_1$ and $\omega_2$ are the fudamental weights of $G_2$ and $\omega_0$ is the fundamental weight of $\operatorname{PSL(2)}$). Then there exists an open
subset of $X$ that is a product of $R_u(Q^-)$ and the irreducible
horospherical $L(Q)$-module of highest weight
$\omega_2-2\omega_1-2\omega_0$.
\end{proof}

\begin{lem}
The two varieties $\mathbf{X_1}$ and $\mathbf{X_2}$ are not homogeneous.
\end{lem}
In fact, this lemma is a corollary of the following lemma, using
the same arguments as in the proof of Lemma \ref{fixedorbit}.
\begin{lem}
The spaces of global sections of
the normal sheafs of the closed orbits $Z_1$ and $Z_2$ in $\mathbf{X_1}$ and $\mathbf{X_2}$ respectively are both trivial.
\end{lem}
\begin{proof}
By the local structure of $\mathbf{X_1}$ and $\mathbf{X_2}$ given in the latter two proofs, one can
remark that the total spaces of the normal sheafs are the
vector bundles $G\times^Q
V(\omega_1-\omega_3)$ and $G\times^Q V(\omega_2-2\omega_1-2\omega_0)$ respectively (and with the respective notations).
We conclude by the Borel-Weil Theorem using the same method as in the
proof of Lemma \ref{fixedorbit}.
\end{proof}

Now, we are going to prove that, in all other cases, a variety
satisfying (*) is homogeneous.

\subsubsection{$G$ simple}

\begin{prop}
In Case (a) of Lemma \ref{touslescas}, $G=\operatorname{PSL(5)}$ and
$X=\Pbb(\bigwedge^2\Cbb^5)$.
\end{prop}

\begin{proof}
Let $G=\operatorname{PSL(5)}$. Let $(e_1,\ldots,e_5)$ be a basis of
$\Cbb^5$. Let $x:=[e_1\wedge e_2+e_3\wedge
e_4]\in\Pbb(\bigwedge^2\Cbb^5)$. Let $H:=\operatorname{Stab}_Gx$.
Then $H$ necessarily stabilizes the subspace $V$ generated by $e_1$,
$e_2$, $e_3$ and $e_4$, so that $H\subset
P:=\operatorname{Stab}_GV$. One can also easily check that
$R(P)\subset H$. Moreover $P/R(P)$, that is isomorphic to
$\operatorname{PSL(4)}$, acts on $\Pbb(\bigwedge^2\Cbb^4)$ with two
orbits as in Case 1b with $n=6$.

We have proved that the orbit of $x$ is isomorphic to the
homogeneous space of Case (a) of Lemma \ref{touslescas} and Case 1b
of Lemma \ref{3cas}. Now let us remark, to conclude, that
$$\dim(G/H)=\dim(G/P)+\dim(P/H)=9.$$
\end{proof}


\begin{prop}\label{casok}
In Case (b) of Lemma \ref{touslescas},
$X=\operatorname{Gr}_q(i+1,n+1)$.
\end{prop}

\begin{proof}
First, suppose that $n=2m\geq 6$. Let
$(e_0,e_1,\ldots,e_m,e_{-1},\ldots,e_{-m})$ be a basis of
$\Cbb^{n+1}$ such that $\langle e_0,e_0\rangle=-2$,  $\langle
e_k,e_l\rangle=1$ if $k=-l\neq 0$ and 0 if $k\neq -l$. Let
$G\simeq\operatorname{SO(n)}$ be the subgroup of
$\operatorname{SO(n+1)}$ stabilizing $e_0$. Let $1\leq i\leq m-2$
and $$\begin{array}{c} x:=[e_1\wedge\cdots\wedge
e_i+e_1\wedge\cdots\wedge e_{i}\wedge
(e_{i+1}+e_{-(i+1)})]\in\Pbb(\bigwedge^i\Cbb^n\oplus\bigwedge^{i+1}\Cbb^n).\end{array}$$
Let $H:=\operatorname{Stab}_Gx$. Then $H$ fixes
$[e_1\wedge\cdots\wedge e_i]$ in $\Pbb(\bigwedge^i\Cbb^n)$ so that
$H\subset P:=P(\omega_{\alpha_i})$. Moreover one can check that
$R(P)\subset H$ and  that the subgroup of $$P/R(P)\simeq
\operatorname{PGL(i)}\times\operatorname{SO(n-2i)}$$ that fixes $x$
is isomorphic to the quotient of
$(\operatorname{PGL(i)}\times\operatorname{SO(n-2i-1)})$ by its
center. Here, $\operatorname{SO(n-2i-1)}$ is the subgroup of
$\operatorname{SO(n-2i)}$ that fixes $e_{i+1}+e_{-(i+1)}$. We have
proved that $G/H$ is an homogeneous space satisfying Case (b) of
Lemma \ref{touslescas} and Case 1a of Lemma \ref{3cas}.

Now let us consider the isomorphism
 $$\begin{array}{ccc}
\bigwedge^i\Cbb^{n}\oplus\bigwedge^{i+1}\Cbb^{n} & \lra & \bigwedge^{i+1}\Cbb^{n+1}\\
x+y & \lmt & x\wedge e_0 +y
\end{array}$$
such that $x$ is map to $[e_1\wedge\cdots\wedge e_i\wedge
(e_{i+1}+e_{-(i+1)}+e_0)]$. Remark that the vector space generated
by $e_1,\ldots,e_i$ and $e_{i+1}+e_{-(i+1)}+e_0$ is isotropic. Then,
by this isomorphism, the closure $X$ of $G.x$ in
$\Pbb(\bigwedge^i\Cbb^n\oplus\bigwedge^{i+1}\Cbb^n)$ is a subvariety
of $\operatorname{Gr}_q(i+1,n+1)$. We conclude saying that $X$ and
$\operatorname{Gr}_q(i+1,n+1)$ have the same dimension.


When $n=2m+1$, the proof is very similar. Indeed, choose
$(e_1,\ldots,e_{m+1},e_{-1},\ldots,e_{-(m+1)})$ be a basis of
$\Cbb^{n+1}$ such that  $\langle e_k,e_l\rangle=1$ if $k=-l\neq 0$
and 0 if $k\neq -l$. Then $G$ is the subgroup of
$\operatorname{SO(n+2)}$ stabilizing $e_0:=e_{m+1}-e_{-(m+1)}$. Then
the remained part of the proof is the same.

\end{proof}

\begin{prop}\label{casc}
In Case (c) of Lemma \ref{touslescas}, $X$ is the homogeneous
variety $F_4/P(\omega_1)$.
\end{prop}

\begin{proof}
Let $G=\operatorname{Spin(9)}\subset F_4$. Denote by $B''$ a Borel
subgroup of $F_4$ and by $B$ the Borel subgroup $B''\cap G$ of $G$.
Let $P:=P(\omega_4)$ and $Q:=P(\omega_2)$ be the corresponding
parabolic subgroups of $G$ containing $B$. Let us denote by
$\omega_i''$ the fundamental weights corresponding to $(F_4,B'')$
and by $\omega_i'$ the fundamental weights corresponding to
$(P/R(P),B/R(P))$ (to fix the numerotation, let us choose that
$P/R(P)$ is of type $D_3$). Consider the homogeneous variety
$F_4/P(\omega_1'')$ as the closed orbit of the projective space
$\Pbb(V(\omega_1''))=\Pbb(V(\omega_4)\oplus V(\omega_2)).$

First, let us prove that $F_4/P(\omega_1'')$ does not contain the
closed orbit $G/P$ of $\Pbb(V(\omega_4))$. Suppose the contrary,
then $P$ is a subgroup of a parabolic group $P''$ of $F_4$
conjugated to $P(\omega_1'')$. In particular, a Levi subgroup of $P$
is included in a Levi subgroup of $P''$. But Levi subgroups of $P$
are semi-direct products of $\Cbb^*$ and a simple group of type
$A_3$, whereas Levi subgroups of $P''$ are semi-direct products of
$\Cbb^*$ and a simple group of type $C_3$. So $P$ cannot be a
subgroup of $P''$, that gives a contradiction.

Remark now that $F_4/P(\omega_1'')$ contains necessarily $G/Q$, the
second closed orbit of \linebreak[4] $\Pbb(V(\omega_4)\oplus
V(\omega_2))$.

Let us consider the rational map $\phi$ from $F_4/P(\omega_1'')$ to
$\Pbb(V(\omega_4))$ defined by projection. It is $G$-equivariant and
its image contains the closed orbit $G/P$ of $\Pbb(V(\omega_4))$.
Let us denote by $X_0'$ the fiber of $\phi$ over the point
$[v_{\omega_4}]$ of $G/P$.  It is stable under the action of $P$ and
its dimension is at most 5 because $\dim(F_4/P(\omega_1''))=15$ and
$\dim(G/P)=10$. The closure of the bundle $G\times^PX_0'\simeq
G.X_0'$ in $\Pbb(V(\omega_4)\oplus V(\omega_2))$ contains
necessarily the closed orbit $G/Q$. Then the closure $X'$ of $X_0'$
contains $P/(P\cap Q)$ (in other words, $[v_{\omega_2}]\in X')$.
Remark that $\dim(P/(P\cap Q))=4$ so that $\dim(X_0')=5$.

Let us decompose $V(\omega_2)$ into a direct sum of $P$-modules
$V'\oplus V''$ where $V'$ is the $P$-module generated by
$v_{\omega_2}$. Then, by the latter paragraph, there exists an
element of $X_0'$ of the form $[v_{\omega_4}+v'+v'']$ where $v'\in
V'\backslash \{0\}$  and $v''\in V''$. Let us remark that
$\omega_2=\omega_4+\omega'_1$, then it is easy to check that there
is no weight, smaller than $\omega_2$, that is the sum of $\omega_4$
and a dominant weight of $(P/R(P),B/R(P))$. It means that the center
of $L(P)$ acts with the same weight on $\Cbb.v_{\omega_4}$ and $V'$
but with a different weight on $V''$. So one can deduce that
$[v_{\omega_4}+v']$ is in $X'$.

Remark now that $V'$ is the $P$-module $\Lambda^2\Cbb^7$, and that
the projective space $\Pbb(\Cbb.v_{\omega_4}\oplus V')$ has 4 types
of $P$-stable subvarieties: the point $[1,0]$, the quadric $P/(P\cap
Q)$ in $\Pbb(V')$, one quadric of dimension 5 with the two latter
closed orbits and infinitely many quadrics of dimension 5 with only
one closed orbit $P/(P\cap Q)$. The first one and the third one
cannot be in $X'$ because $G/P$ is not included in
$F_4/P(\omega_1'')$. It is easy to see that the second one cannot
also be in $X'$ (because $P/(P\cap Q)$ is not in $X_0'$). Then, we
conclude that $X'$ is a quadric as in Case 1a of Lemma \ref{3cas}
with $n=6$. And $X:=\overline{G.X'}\subset
F_4/P(\omega_1'')\subset\Pbb(V(\omega_4)\oplus V(\omega_2))$ is a
variety with an open orbit satisfying Case (c) of Lemma
\ref{touslescas} and Case 1a of Lemma \ref{3cas}. Since the
dimension of $X$ and $F_4/P(\omega_1'')$ are the same, we have
necessarily  $X=F_4/P(\omega_1'')$.
\end{proof}

\begin{prop}
In Case (d) of Lemma \ref{touslescas}, $X$ is a quadric of
dimension~14 in the projectivization of the spinorial
$\operatorname{Spin(9)}$-module.
\end{prop}

\begin{proof}
In fact $X$ is the quadric $g(x)=0$ of \cite[Prop.5]{Ig70}.
\end{proof}

\begin{prop}\label{newcase}
In Case (e) of Lemma \ref{touslescas}, there exists no variety
satisfying (*).
\end{prop}

\begin{proof}
In that case $G=\operatorname{Spin(7)}$ (or $\operatorname{SO(7)}$),
and $P=P(\omega_2)$. Then $P/R(P)$ is isomorphic to
$\operatorname{PSL(2)}\times\operatorname{SL(2)}$ (or
$\operatorname{PSL(2)}\times\operatorname{SL(2)}$).

Suppose there exists a variety $X$ satisfying (*) in Case (e) of
\ref{touslescas}. Then, by the preceding paragraph, the associated
variety $X'$ must be the one of Case 1b of Lemma \ref{3cas} with
$n=4$.


Let us now check that $X$ is not smooth, in order to obtain a
contradiction. By Remark \ref{structurelocaleX'}, there exists an
affine open subset of $X'$ that is a product of $R_u((P\cap Q)^-)$
and a line where $P\cap Q$ acts with weight
$2\omega_2-2\omega_1-2\omega_3$ (it is Lemma \ref{smooth} in this
particular case). Then there exists an open subset of $X$ that is a
product of $R_u(Q^-)$ and a cone on the flag variety in the
irreducible $L(Q)$-module of highest weight
$2\omega_2-2\omega_1-2\omega_3$. This implies that $X$ is not
smooth.


\end{proof}



\begin{prop}
In Case (f) of Lemma \ref{touslescas}, $X$ is the grassmannian
$\operatorname{Gr}(3,2n)$.
\end{prop}

\begin{proof}
Let $G=\operatorname{Sp(2n)}$. First, let us consider the
grassmannian $\operatorname{Gr}(3,2n)$ under the action of $G$. It
has clearly two orbits: the open one $X_0$ consisting of
non-isotropic elements of $\operatorname{Gr}(3,2n)$ and the closed
one isomorphic to
$G/P(\omega_3)\simeq\operatorname{Gr}_\omega(3,2n)$. We just have to
compute the stabilizer of a point of $X_0$ and to check that it is
the good one. Let $(e_1,\ldots,e_n,e_{-n},\ldots,e_{-1})$ be a basis
of $\Cbb^{2n}$ such that $\langle e_k, e_l \rangle =1$ if $k=-l>0$,
$-1$ if $k=-l<0$ and $0$ if $k\neq -l$. Denote by $E$ the subspace
of $\Cbb^{2n}$ generated by $e_1$, $e_2$ and $e_{-2}$. Let
$H:=\operatorname{Stab}_GE$. Remark that $H\subset
P:=\operatorname{Stab}_G[e_1]$ because the line generated by $e_1$
is the only line of $E\cap E^\perp$. Then one can check that
$R(P)\subset H$ and  that $P/H\simeq \operatorname{Sp(
 2n-2)}/(\operatorname{Sp(2)}\times\operatorname{Sp(2n-4})).$
\end{proof}

\begin{prop}\label{newcasg}
In Case (g) of Lemma \ref{touslescas}, there exists no variety
satisfying (*).
\end{prop}

\begin{proof}
Suppose that there is a variety $X$ satisfying (*) in Case (g) of
Lemma \ref{touslescas}. With the same method as in the proof of
Proposition \ref{newcase}, we prove that
there exists an open subset of $X$ that is a product of $R_u(Q^-)$
and a cone on the flag variety in the irreducible $L(Q)$-module of
highest weight $3\omega_2-2\omega_1-2\omega_3$. This implies that
$X$ is not smooth.
\end{proof}

\begin{prop}
In Case (i) of Lemma \ref{touslescas}, $X$ is the homogeneous
variety $E_6/P(\omega_2)$.
\end{prop}

\begin{proof}
We use exactly the same method as in the proof of Proposition
\ref{casc}, so that we only give here a sketch of the proof.

Let $G=F_4\subset E_6$, $P:=P(\omega_4)$ and $Q:=P(\omega_1)$. We
consider $E_6/P(\omega_2'')$ in the projective space
$\Pbb(V(\omega_2''))=\Pbb(V(\omega_4)\oplus V(\omega_1)).$

First,  $E_6/P(\omega_2'')$ does not contain the closed orbit $G/P$
of $\Pbb(V(\omega_4))$, because a Levi subgroup of $P$ cannot be
included in a Levi subgroup of a group conjugated to $P(\omega_2')$.

Moreover $E_6/P(\omega_2'')$ contains necessarily $G/Q$.

Now we consider the rational map $\phi$ from $E_6/P(\omega_2'')$ to
$\Pbb(V(\omega_4))$ defined by projection. It is $G$-equivariant and
its image contains the closed orbit $G/P$ of $\Pbb(V(\omega_4))$.
Let denote by $X_0'$ the fiber of $\phi$ over the point
$[v_{\omega_4}]$ of $F_4/P$.  It is stable under the action of $P$
and its dimension is at most 6 because $\dim(E_6/P(\omega_2''))=21$
and $\dim(G/P)=15$.

Then one can prove that the closure $X'$ of $X_0'$ contains
$P/(P\cap Q)$, and deduce that $X'$ contains an element of the form
$[v_{\omega_4}+v']$, where $v'$ is a non-zero element of the
sub-$P$-module $V'$ generated by $v_{\omega_1}$.

Remark now that $V'$ is the $P$-module $\Cbb^7$, and that the
projective space $\Pbb(\Cbb.v_{\omega_4}\oplus V')$ has 4 types of
$P$-stable subvarieties: the point $[1,0]$, the quadric $P/(P\cap
Q)$ in $\Pbb(V')$, one quadric of dimension 6 with the two latter
closed orbits and infinitely many quadrics of dimension 6 with only
one closed orbit $P/(P\cap Q)$.

Then we conclude from the latter paragraph that $X'$ is the
6-dimensional quadric of Case 1a of Lemma \ref{3cas} with $n=7$, so
that $E_6/P(\omega_2'')\subset\Pbb(V(\omega_4)\oplus V(\omega_1))$
is the variety $X$ satisfying (*) in Case (i) of Lemma
\ref{touslescas}.

\end{proof}

\begin{prop}\label{cash}
In Case (j) of Lemma \ref{touslescas}, there exists no variety
satisfying (*).
\end{prop}

\begin{proof}
Suppose that there is a variety $X$ satisfying (*) in Case (j) of
Lemma \ref{touslescas}. With the same method as in the proof of
Proposition \ref{newcase}, we prove that
there exists an open subset of $X$ that is a product of $R_u(Q^-)$
and a cone on the flag variety in the irreducible $L(Q)$-module of
highest weight $3\omega_4-2\omega_3$. This implies that $X$ is not
smooth.
\end{proof}

\subsubsection{$G$ not simple}

\begin{prop}
In Case (a') of Lemma \ref{casnotsimple},
$G=\operatorname{PSL(n+1)}\times\operatorname{PSL(2)}$ and $X$ is
the projective space $\Pbb(\Cbb^{n+1}\otimes\Cbb^2)$.
\end{prop}

\begin{proof}
Let $G=\operatorname{PSL(n+1)}\times\operatorname{PSL(2)}$ acting on
$\Pbb(\Cbb^{n+1}\otimes\Cbb^2)$.

 Let $(e_1,\ldots,e_{n+1})$ and $(f_1,f_2)$ be respectively some basis of
$\Cbb^{n+1}$ and $\Cbb^2$. Let $x:=[e_1\otimes f_1 +e_2\otimes
f_2]\in\Pbb(\Cbb^{n+1}\otimes\Cbb^2)$ and
$H:=\operatorname{Stab}_Gx$.

Then $H$ is clearly included in the maximal parabolic subgroup $P$
of $G$ that stabilizes the plane generated by $e_1$ and $e_2$.
Moreover, one can check that $R(P)\subset H$. Remark that $P/H$ is
isomorphic to
$(\operatorname{PSL(2)}\times\operatorname{PSL(2)})/\operatorname{PSL(2)}$
where $\operatorname{PSL(2)}$ is included in
$\operatorname{PSL(2)}\times\operatorname{PSL(2)}$ as in proof of
Proposition \ref{casf'}.

We have proved that $G/H$ is the homogeneous space of Case (a') of
Lemma \ref{casnotsimple}. We complete the proof saying that the
dimension of $G/H$ is $2n+1$.
\end{proof}

\begin{prop}
In Case (b') of Lemma \ref{casnotsimple}, $X$ is the spinorial
variety $\operatorname{Gr}^+_q(n+2,2n+4)$.
\end{prop}

\begin{proof}
Let $(e_1,\ldots,e_{n+2},e_{-1},\ldots,e_{-(n+2)})$ be a basis of
$\Cbb^{2n+4}$ such that $\langle e_k,e_l\rangle=1$ if $k=-l$ and $0$
if $k\neq -l$. Let
$G\simeq\operatorname{SO(2n+1)}\times\operatorname{SO(3)}$ be the
subgroup of $\operatorname{SO(2n+4)}$ that stabilizes the two
subspaces $V_1$ and $V_2$ of $\Cbb^{2n+4}$ respectively generated by
$$e_1,\ldots,e_n,e_1,\ldots,e_{-n},e_{n+1}+e_{-(n+1)}\mbox{ and }e_{n+1}-e_{-(n+1)},e_{n+2},e_{-(n+2)}.$$

Let $V\subset \Cbb^{2n+4}$ be the $(n+2)$-dimensional subspace
generated by
$e_1,\ldots,e_{n-1},e_{n+1},e_n+e_{n+2},e_{-n}+e_{-(n+2)}$ and
$H:=\operatorname{Stab}_GV$. Then $H$ is included in the maximal
parabolic subgroup $P$ of $G$  that stabilizes the
$(n-1)$-dimensional subspace $V\cap V_1$. One can check also that
$R(P)\subset H$. Now, $P/H$ is isomorphic to
$(\operatorname{SO(3)}\times\operatorname{SO(3)})/H'$ where $H'$ is
the stabilizer of the 3-dimensional subspace generated by
$e_{-(n+1)},e_n+e_{n+2},e_{-n}-e_{-(n+2)}$ in the group
$\operatorname{SO(3)}\times\operatorname{SO(3)}\subset G$ acting on
the subspace $W$ generated by
$e_n,e_{n+1},e_{n+2},e_{-n},e_{-(n+1)}, \linebreak[4] e_{-(n+2)}$.
Remark that $2e_{-(n+1)}=(e_{n+1}+e_{-(n+1)})-(e_{n+1}-e_{-(n+1)})$,
so that $H'$ is the group $\operatorname{SO(3)}$ diagonally embedded
in $\operatorname{SO(3)}\times\operatorname{SO(3)}$.

We complete the proof saying that $G/H$ and
$\operatorname{Gr}^+_q(n+2,2n+4)$ have the same dimension.

\end{proof}

\begin{prop}
In Case (c') of Lemma \ref{casnotsimple}, $X$ is the symplectic
grassmannian $\operatorname{Gr}_\omega(n+1,2n+2)$.
\end{prop}

\begin{proof}
Let $V$ and $V'$ be the fundamental $\operatorname{Sp(2n)}$-modules
respectively of weight $\omega_n$ and $\omega_{n-1}$. Let
$G=(\operatorname{Sp(2n)}\times\operatorname{SL(2)})$ acting on
$\Pbb(V'\oplus (V\otimes\Cbb^2))$.

Let $(e_1,\ldots,e_n,e_{-1},\ldots,e_{-n})$ and
$(e_{n+1},e_{-(n+1)})$ be respectively some basis of $\Cbb^{2n}$ and
$\Cbb^2$, such that $\langle e_k,e_l\rangle=1$ if $k=-l>0$, $\langle
e_k,e_l\rangle=-1$ if $k=-l<0$ and $0$ if $k\neq -l$.  Remark that
$\langle\,,\,\rangle$ can be naturrally defined on
$\Cbb^{2n}\oplus\Cbb^2$ so that $G$ is a subgroup of
$\operatorname{Sp(2n+2)}$. Note also that
$V\subset\bigwedge^n\Cbb^{2n}$ and
$V'\subset\bigwedge^{n-1}\Cbb^{2n}$. Let $$x:=[e_1\wedge\ldots\wedge
e_{n-1}+(e_1\wedge\ldots\wedge e_n)\otimes e_{n+1}
-(e_1\wedge\ldots\wedge e_{n-1}\wedge e_{-n})\otimes e_{-(n+1)}]$$
and $H:=\operatorname{Stab}_Gx$.

Then $H$ is clearly included in the maximal parabolic subgroup $P$
of $G$ that stabilizes the subspace generated by $e_1,\ldots,
e_{n-1}$. Moreover, one can check that $R(P)\subset H$. Remark that
$P/R(P)=(\operatorname{Sp(2n-2)}\times\operatorname{SL(2)}\times\operatorname{SL(2)})$
so that $P/H$ is isomorphic to
$(\operatorname{SL(2)}\times\operatorname{SL(2)})/\operatorname{SL(2)}$
where $\operatorname{SL(2)}$ is included in
$\operatorname{SL(2)}\times\operatorname{SL(2)}$ as the same way as
in the proof of Proposition \ref{P5}.

Now let us remark that the $G$-module $V'\oplus V\otimes\Cbb^2$ is
isomorphic to the fundamental $\operatorname{Sp(2n+2)}$-module $V''$
of weight $\omega_{n+1}$. Indeed, we have the following
$G$-equivariant isomorphism
$$\begin{array}{ccc}
V'\oplus V\otimes\Cbb^2 & \lra & V''\subset\bigwedge^{n+1}\Cbb^{2n+2}\\
y+z\otimes w & \lmt & y\wedge e_{n+1}\wedge e_{-(n+1)}-y\wedge
e_1\wedge e_{-1}-\cdots -y\wedge e_n\wedge e_{-n}+z\wedge w.
\end{array}$$

Moreover, this isomorphism maps $x$ to
$$\displaylines{e_1\wedge\ldots\wedge e_{n+1}\wedge
e_{-(n+1)}-e_1\wedge\ldots\wedge e_n\wedge
e_{-n}+e_1\wedge\ldots\wedge e_n\wedge e_{n+1}-e_1\wedge\ldots\wedge
e_{n-1}\wedge e_{-n}\wedge e_{-(n+1)} \cr =e_1\wedge\ldots\wedge
e_{n-1}\wedge(e_n-e_{-(n+1)})\wedge(e_{n+1}-e_{-n}).}$$ Remark now
that the vector space generated by $e_1,\ldots,e_{n-1}$,
$e_n-e_{-(n+1)}$ and $e_{n+1}-e_{-n}$ is isotropic so that $G/H$ is
contained in $\operatorname{Gr}_\omega(n+1,2n+2)$.

We complete the proof with an argument of dimension.
\end{proof}

\begin{prop}
In Case (d') of Lemma \ref{casnotsimple},
$G=\operatorname{Sp(2n)}/\{\pm 1\}\times\operatorname{PSL(2)}$ and
$X$ is a quadric in $\Pbb(\Cbb^{2n}\otimes\Cbb^2)$.
\end{prop}

\begin{proof}
Let $G=\operatorname{Sp(2n)}/\{\pm 1\}\times\operatorname{PSL(2)}$,
acting on $\Pbb(\Cbb^{2n}\otimes\Cbb^2)$.

Let $(e_1,\ldots,e_n,e_{-1},\ldots,e_{-n})$ and $(f_1,f_2)$ be
respectively some basis of $\Cbb^{2n}$ and $\Cbb^2$, such that
$\langle e_k,e_l\rangle=1$ if $k=-l>0$, $\langle e_k,e_l\rangle=-1$
if $k=-l<0$ and $0$ if $k\neq -l$. Let $x:=[e_1\otimes f_1
+e_2\otimes f_2]\in\Pbb(\Cbb^{2n}\otimes\Cbb^2)$ and
$H:=\operatorname{Stab}_Gx$.

Then $H$ is clearly included in the parabolic subgroup $P$ of $G$
that stabilizes the isotropic plane generated by $e_1$ and $e_2$.
Moreover, one can check that $R(P)\subset H$. Then remark that $P/H$
is isomorphic to
$(\operatorname{PSL(2)}\times\operatorname{PSL(2)})/\operatorname{PSL(2)}$
where $\operatorname{PSL(2)}$ is included in
$\operatorname{PSL(2)}\times\operatorname{PSL(2)}$ as in the proof
of Proposition \ref{P5}.

Now, let us define a quadratic form $q$  on $\Cbb^{2n}\otimes\Cbb^2$
by $$\forall x_1,\,x_2\in \Cbb^{2n},\,q(x_1\otimes f_1+x_2\otimes
f_2)=\langle x_1,x_2\rangle.$$ One can check that $q$ is invariant
under the action of $G$. Moreover $q(e_1\otimes f_1 +e_2\otimes
f_2)=0$ so that $G/H$ is contained in the quadric
$Q^{4n-2}:=\{[y]\in\Pbb^{4n-1}\mid q(y)=0\}$. To conclude, note that
$\dim(G/H)=4n-2$.
\end{proof}

\begin{prop}\label{case'}
In Case (e') of Lemma \ref{casnotsimple}, there exists no variety
satisfying (*).
\end{prop}

\begin{proof}
Let us denote by $\omega_0$ the fundamental weight of
$\operatorname{PSL(2)}$. Suppose that there is a variety $X$
satisfying (*) in Case (e') of Lemma \ref{casnotsimple}. With the
same method as in the proof of Proposition \ref{newcase}, we prove
that
there exists an open subset of $X$ that is a product of $R_u(Q^-)$
and a cone on the flag variety in the irreducible $L(Q)$-module of
highest weight $3\omega_1-2\omega_2-2\omega_0$. This implies that
$X$ is not smooth.
\end{proof}

\subsection{When $R(P)\not\subset H$}\label{casiiiv}

When $H$ contains a Levi subgroup of $P$, we have the following
result.

\begin{lem}[Th.2.1 of \cite{Br88}]\label{4cas}

Suppose that $H$ contains a Levi subgroup of the maximal parabolic
subgroup $P$.

Then $(G,P,H)$ is one of the following:

(i) $(\operatorname{PSL(m+1)},P(\omega_1),\operatorname{GL(m)})$.

(ii) $G=\operatorname{SO(2m+1)}$; $P=P(\omega_m)$ is the stabilizer
of an isotropic $m$-dimensional subspace $E$ of $\Cbb^{2m+1}$;  $H$
is the stabilizer of $E$ and a non-isotropic vector orthogonal to
$E$.

(iii) $G=\operatorname{Sp(2m)}/\{\pm 1\}$; $P$ is the stabilizer of
a line $l$ in $\Cbb^{2m}$; $H$ is the stabilizer of $l$ and a
non-isotropic plane containing $l$.

(iv) $G=G_2$; $P$ is the stabilizer of a line $l$ in the
7-dimensional simple $G$-module;  $H$ is the stabilizer of a line
$l'$ in the 14-dimensional simple $G$-module such that $H$ contains
a maximal torus $T$ of $G$ and $T$ acts with the same weight
$\alpha$ in $l$ and $l'$; $\alpha$ is a short root of $(G,T)$.
\end{lem}

In Case (i), $\tilde{X}=\Pbb^n\times (\Pbb^n)^*$, and in Case (iii),
it is easy to check that $\tilde{X}$ is the partial flag variety
$\operatorname{SL(2m)}/(P(\omega_1)\cap P(\omega_2))$. So these two
cases cannot give us a variety satisfying (*). And in the next two
propositions we are going to prove that $\tilde{X}$ is homogeneous
in Case (ii) and that $X$ is homogeneous in Case (iv).

\begin{prop}
Let $\tilde{X}$ be the $\operatorname{SO(2n+1)}$-variety defined by
$$\tilde{X}:=\{(l,V)\in \operatorname{Gr}(1,2n+1)\times \operatorname{Gr}_q(n,2n+1)\mid l\subset V^\perp\}$$ where $V^\perp$ is the $q$-orthogonal subspace of $V$.

Then $\tilde{X}\simeq \operatorname{SO(2n+2)}/(P(\omega_1)\cap
P(\omega_{n+1}))$.
\end{prop}

\begin{rem}
The variety $\tilde{X}$ is the two-orbits variety with an open orbit
isomorphic to the homogeneous space of Case (ii) of Lemma \ref{4cas}
and with closed orbit of codimension~1.
\end{rem}

\begin{proof}
We have to prove that $$\tilde{X}\simeq \tilde{Y}:= \{(l',V')\in
\operatorname{Gr}_q(1,2n+2)\times
\operatorname{Gr}_q^+(n+1,2n+2)\mid l'\subset V'\}.$$ Let us
decompose $\Cbb^{2n+2}$ in an orthogonal sum $\Cbb^{2n+1}\oplus L$
such that the restriction of $q$ to $\Cbb^{2n+1}$ is of maximal
index. Denote by $\pi$ the orthogonal projection
$\Cbb^{2n+2}\lra\Cbb^{2n+1}$.

Let $(l',V')\in \tilde{Y}$. Then $\pi(l')$ is a line because $L$ is
not isotropic. Also $V'\cap\Cbb^{2n+1}$ is an isotropic subspace of
$\Cbb^{2n+1}$ of dimension $n$. Moreover $\pi(l')\subset
(V'\cap\Cbb^{2n+1})^\perp$. Indeed let $v$ a non-zero element of
$\pi(l')$. One can write $v=v_1+v_2$ with $v_1\in l'$ and $v_2\in
L$. Then $v_1\in V'^\perp$ and $v_2\in (\Cbb^{2n+1})^\perp$, so that
$v\in(V'\cap\Cbb^{2n+1})^\perp$.

By the latter paragraph, one can define a
$\operatorname{SO(2n+1)}$-equivariant morphism
$$\begin{array}{cccc}
\phi : & \tilde{Y} & \lra & \tilde{X} \\
 & (l',V') & \lmt & (\pi(l'), V'\cap\Cbb^{2n+1})
\end{array}$$

Let us show now that $\phi$ is an isomorphism.

Let $(l,V)\in \tilde{X}$. Then there exists a unique $V'\in
\operatorname{Gr}_q^+(n+1,2n+2)$ such that $V\subset V'$ and
$V=V'\cap\Cbb^{2n+1}$.

By hypothesis $(l\oplus L)\subset V^\perp\subset\Cbb^{2n+2}$. But
$V'$ is also a subspace of $V^\perp$ which is of dimension $n+2$.
Moreover $l\oplus L$ is not isotropic so it is not included in $V'$.
So $l':=(l\oplus L)\cap V'$ is an isotopic line in $V'$ and
$\pi(l')=l$. Let us note also that a such $l'$ is unique (because it
must be included in $V'$ and in $l\oplus L$).
\end{proof}

\begin{prop}\label{casiv}
The variety $\operatorname{Gr}_q(2,7)$ has two orbits under the
action of $G_2$. It is the two-orbits variety satisfying (*) in Case
(iv) of Lemma \ref{4cas}.
\end{prop}

\begin{proof}
Let $\Obb$ be the set of octonions on $\Cbb$.
One can define this 8-dimensional (non-associative) algebra with a basis $(1,e_1,\ldots,e_7)$ and the following multiplication table ($1$ being the identity):
\begin{center}
\begin{tabular}{|c||c|c|c|c|c|c|c|}
\hline
 & $e_1$ & $e_2$ & $e_3$ & $e_4$ & $e_5$ & $e_6$ & $e_7$\\
\hline
\hline
 $e_1$ & $ -1$ & $e_4$ & $e_7$ & $-e_2$ & $e_6$ & $-e_5$ & $-e_3$\\
\hline
 $e_2$ & $-e_4$ & $-1$ & $e_5$ & $e_1$ & $-e_3$ & $e_7$ & $-e_6$\\
\hline
 $e_3$ & $-e_7$ & $-e_5$ & $-1$ & $e_6$ & $e_2$ & $-e_4$ & $e_1$\\
\hline
 $e_4$ & $e_2$ & $-e_1$ & $-e_6$ & $-1$ & $e_7$ & $e_3$ & $-e_5$\\
\hline
$e_5$  & $-e_6$ & $e_3$ & $-e_2$ & $-e_7$ & $-1$ & $e_1$ & $e_4$\\
\hline
$e_6$ & $e_5$ & $-e_7$ & $e_4$ & $-e_3$ & $-e_1$ & $-1$ & $e_2$\\
\hline
$e_7$ & $e_3$ & $e_6$ & $-e_1$ & $e_5$ & $-e_4$ & $-e_2$ & $-1$\\
\hline
\end{tabular}
\end{center}

Let us denote by $q$ the norm defined by
$q(x_0+\sum_{i=1}^7x_ke_k)=\sum_{k=0}^7 x_k^2$ and $\Im(\Obb)$ the
7-dimensional subspace generated by $e_1,\ldots,e_7$. Let $SO(7)$ be
the special orthogonal group defined by $q$ and acting on
$\Im(\Obb)$. Let $G$ be the group of automorphism of $\Obb$. Remark
that, since an automorphism fixes the identity and preserve the
norm, $G$ is a subgroup of $SO(7)$. Moreover it is well-known that
$G$ is of type $G_2$.

Let $\operatorname{Gr}_q(2,7)$ the set of isotropic planes in
$\Im(\Obb)$. First, remark that $\operatorname{Gr}_q(2,7)$ has a
closed orbit $Z$ isomorphic to $G/P(\omega_2)$.  In fact $Z=\{V\in
\operatorname{Gr}_q(2,7)\mid \forall z,z'\in V,\,zz'=0\}$.

Let us compute the stabilizer of a particular point of
$\operatorname{Gr}_q(2,7)$. For this, let us consider a new basis
$(z_0,z_1,z_2,z_3,z_{-1},z_{-2},z_{-3})$ of $\Im(\Obb)$ having the
following multiplication table.
\begin{center}
\begin{tabular}{|c||c|c|c|c|c|c|c|}
\hline
& $z_0$ & $z_1$ & $z_2$ & $z_3$ & $z_{-1}$ & $z_{-2}$ & $z_{-3}$ \\
\hline
\hline
 $z_0$ & $ 1$ & $z_1$ & $z_2$ & $-z_3$ & $-z_{-1}$ & $-z_{-2}$ & $z_{-3}$\\
\hline
 $z_1$ & $-z_1$ & 0 & $z_3$ & 0 & $-1-z_0$ & 0 & $-2z_{-2}$\\
\hline
 $z_2$ & $-z_2$ & $-z_3$ & 0 & 0 & 0 & $-1-z_0$ & $2z_{-1}$\\
\hline
 $z_3$ & $z_3$ & 0 & 0 & 0 & $2z_2$ & $-2z_1$ & $-2+2z_0$\\
\hline
$z_{-1}$  & $z_{-1}$ & $-1+z_0$ & 0 & $-2z_2$ & 0 & $z_{-3}$ & 0\\
\hline
$z_{-2}$ & $z_{-2}$ & 0 & $-1+z_0$ & $2z_1$ & $-z_{-3}$ & 0 & 0\\
\hline
$z_{-3}$ & $-z_{-3}$ & $2z_{-2}$ & $-2z_{-1}$ & $-2-2z_0$ & 0 & 0 & 0\\
\hline
\end{tabular}
\end{center}
Take for example, $z_0=\imath e_7$, $z_1=(e_1+\imath e_3)/\sqrt{2}$,
$z_2=(e_2+\imath e_6)/\sqrt{2}$, $z_3=e_4-\imath e_6$,
$z_{-1}=(e_1-\imath e_3)/\sqrt{2}$, $z_{-2}=(e_2-\imath
e_6)/\sqrt{2}$ and $z_{-3}=e_4+\imath e_6$.

Let $E$ be the plane of $\Im(\Obb)$ generated by $z_1$ and $z_2$ and
let $H:=\operatorname{Stab}_GE$. Remark that $H$ contains a maximal
torus $T$ (the diagonal matrices of $G$ in the new basis). Moreover,
$H$ stabilizes the line $l$ generated by $z_3$ (because
$z_1z_2=z_3$). In other words, $H\subset P:=\operatorname{Stab}_Gl$.
Moreover $T$ acts on $l$ with weight a short root $\alpha$ of
$(G,T)$. Let us note that the adjoint $G$-module $V^{14}$ is the
submodule of the adjoint $\operatorname{SO(7)}$-module
$\bigwedge^2\Im(\Obb)$. Indeed, it is the kernel of the following
linear map.
$$\begin{array}{cccc}
\bigwedge^2\Im(\Obb) & \lra & \Im(\Obb) \\
 z\wedge z' & \lmt & \Im(zz')
\end{array}$$
Then, the only line of $V^{14}$ where $T$ acts with weight $\alpha$
is the line $l'$  generated by $z_1\wedge z_2+z_0\wedge z_3$. Let us
prove that $H$ is the stabilizer of $l'$. Let $\phi\in H$. Then
$\phi(z_1\wedge z_2)=D.z_1\wedge z_2$ where $D$ is the determinant
of the restriction to $E$ of $\phi$. Moreover, we clearly have
$\phi(z_3)=D.z_3$ (because $z_3=z_1z_2$). Since $z_0$ fixes each
point of $E$ by left multiplication, $\phi(z_0)$ must be of the form
$z_0+\lambda z_3$ for some $\lambda\in\Cbb$. So $\phi(z_0\wedge
z_3)=D.z_0\wedge z_3$ and $\phi(l')=l'$.

Then we have proved that $G/H$ is isomorphic to the homogeneous
space of Case (iv) of Lemma \ref{4cas}.


We complete the proof saying that $G/H$ and
$\operatorname{Gr}_q(2,7)$ have the same dimension.
\end{proof}


\begin{thebibliography}{111111}


\bibitem{Ah83} D.~N.~Akhiezer, {\sl Equivariant completions of homogeneous algebraic varieties by homogeneous divisors},
Ann. Global  Anal. Geom. {\bf 1} (1983), no. 1, 49-78.

\bibitem{Ak95} D.~N.~Akhiezer, {\sl Lie group actions in complex analysis}, Aspects of Mathematics, E27. Friedr. Vieweg and Sohn, Braunschweig, 1995.


\bibitem{Bo75} N.~Bourbaki, {\sl Groupes et alg\`ebres de Lie}, Chapitres 4,5,6, C.C.L.S., Paris 1975.


\bibitem{Br89} M.~Brion, {\sl Groupe de Picard et nombres caract\'eristiques des vari\'et\'es sph\'eriques}, Duke Math. J. {\bf 58} (1989), no. 2, 397-424.

\bibitem{Br88} M.~Brion, {\sl On spherical varieties of rank one (after D. Akhiezer, A. Huckleberry, D. Snow)},
Group actions and invariant theory (Montreal, PQ, 1988),  31-41, CMS
Conf. Proc., {\bf 10}, Amer. Math. Soc., Providence, RI, 1989.

\bibitem{Br97b} M.~Brion, {\sl Vari\'et\'es sph\'eriques}, lecture notes available at http://www-fourier.ujf-grenoble.fr/~mbrion/spheriques.pdf, 1997.




\bibitem{CF03} S.~Cupit-Foutou, {\sl Classification of two-orbit varieties}, Comment. Math. Helv. {\bf 78} (2003) 245-265.

\bibitem{Ig70} J.~Igusa, {\sl A classification of spinors up to dimension twelve},
Amer. J. Math. {\bf 92} (1970) 997--1028.

\bibitem{Kn91} F.~Knop, {\sl The Luna-Vust Theory of Spherical Embeddings}, Proceedings of the Hyderabad Conference on Algebraic Groups, Manoj-Prakashan, 1991, 225-249.

\bibitem{LV83} D.~Luna and T.~Vust, {\sl Plongements d'espaces homog\`enes}, Comment. Math. Helv. {\bf 58} (1983), 186-245.

\bibitem{Mi05} I.~A.~Mihai {\sl Odd symplectic flag manifolds}, preprint, arXiv: math.AG/0604323.

\bibitem{Pa06} B.~Pasquier, {\sl Vari\'et\'es horosph\'eriques de Fano}, thesis available at
http://tel.archives-ouvertes.fr/tel-00111912.

\bibitem{PV94} V.~Popov and E.~Vinberg, {\sl Invariant Theory}, Encyclopedia of Mathematical Sciences, vol 55, Springer Verlag 1994, 123-278.

\bibitem{Ru07} A.~Ruzzi, {\sl Smooth Projective Symmetric Varieties with Picard Number equal to one}, preprint, arXiv: math.AG/0702340.

\bibitem{Sp98} T.A.~Springer, {\sl Linear Algebraic Groups, Second Edition}, Birkh\"auser, 1998.

\bibitem{Ti06} D.~Timashev, {\sl Homogeneous spaces and equivariant embeddings}, preprint, arXiv: math.AG/0602228.

\end{thebibliography}
\end{document}